\input amssym.tex
%%%%%%%%%%%%%%%%%%%%%%%
\font\twbf=cmbx12

\font\sc=cmcsc10

\def\rond{{\scriptstyle\circ}}

\def\cf{{\it cf.\/}\ }    
\def\ie{{\it i.e.\/}\ }
 
\def\up#1{\raise 1ex\hbox{\sevenrm#1}}

\def \bg {\bigskip \goodbreak}
\def \sn {\nobreak \smallskip}

\def\ref#1&#2&#3&#4&#5\par{\par{\leftskip = 5em {\noindent
\kern-5em\vbox{\hrule height0pt depth0pt width
5em\hbox{\bf[\kern2pt#1\unskip\kern2pt]\enspace}}\kern0pt}
{\sc\ignorespaces#2\unskip},\
{\rm\ignorespaces#3\unskip}\
{\sl\ignorespaces#4\unskip\/}\
{\rm\ignorespaces#5\unskip}\par}}

\def \eps {\varepsilon}

\def \reel{\Bbb{R} }
\def \comp{ \Bbb{C} }
\def \nat{ \Bbb{N}}
 
\def\ent{ \Bbb{Z}}

\def \adh  {\overline}

\def \abs#1{\left\vert#1\right\vert }

\def\N#1{\muskip0=-2mu{\left|\mkern\muskip0\left|
#1\right|\mkern\muskip0\right|}}

\def\NN#1{\muskip0=-2mu{\left|\mkern\muskip0\left|
\mkern\muskip0\left|#1\right|\mkern\muskip0
\right|\mkern\muskip0\right|}}

\def\pe{\Bbb{P}}

\def\wa{{\rm 1\mkern-2mu}\!{\rm I}}
\def\tore{\Bbb{T}}
\def\hot{\widehat\otimes}
\def\wt#1{\widetilde #1}
\def\wcheck#1{\smash{
        \mathop{#1}\limits^{\scriptscriptstyle{\vee}}}}

\def\Otimes_#1^#2{\matrix{{}_{#2}\cr
      \bigotimes \cr {}^{#1}}}
\def\hb{\hfill\break}  \def\hbi{\hb\indent}

\def\wco#1#2{\lower0.5ex\hbox{$\scriptstyle{{#1}\wcheck\otimes
{#2}}$}}

\def\refl#1&#2&#3&#4&#5\par{\par{\leftskip = 7em {\noindent
\kern-7em\vbox{\hrule height0pt depth0pt width
5em\hbox{\bf[\kern2pt#1\unskip\kern2pt]\enspace}}\kern0pt}
{\sc\ignorespaces#2\unskip},\
{\rm\ignorespaces#3\unskip}\
{\sl\ignorespaces#4\unskip\/}\
{\rm\ignorespaces#5\unskip}\par}}

\def\ref#1&#2&#3&#4&#5\par{
\item{[{\bf\ignorespaces#1\unskip}]}
{\sc\ignorespaces#2\unskip},\
{\rm\ignorespaces#3\unskip}\
{\sl\ignorespaces#4\unskip\/}\
{\rm\ignorespaces#5\unskip}\par}

\def\th{\theta}

\def\liminf_#1{\mathop{\underline{\hbox{\rm lim}}}\limits_{#1}}
\def\limsup_#1{\mathop{\overline{\hbox{\rm lim}}}\limits_{#1}}

%%%%%%%%%%%%%%%%%%%%%%%
\magnification=1200

\centerline{\twbf On  The Interpolation  of Injective
or  Projective }\bg
\centerline{\twbf Tensor Products  of Banach Spaces }
\bg
\bg
\centerline{{\sc Omran Kouba}}\par
\centerline{Department of Mathematics}\par
\centerline{{\sl Higher Institute for Applied Sciences and Technology}}\par
\centerline{P.O. Box 31983, Damascus, Syria.}\par
\centerline{{\it E-mail} : omran\_kouba@hiast.edu.sy}\par
\bg
\baselineskip=13pt
\bg
{\bf Abstract :} 
We prove a general  result  on  the  factorization of  
matrix-valued analytic functions. 
We deduce that  if $(E_0,E_1)$ and $(F_0,F_1)$
are interpolation 
 pairs  with dense  intersections, then under some 
conditions  on the  spaces $E_0$, $E_1$, $F_0$  and $
F_1$, we have
$$
[E_0\hot F_0 , E_1\hot F_1 ]_\th =
[E_0 ,E_1]_\th\hot [F_0 ,F_1]_\th,
\qquad 0 < \th < 1.$$
\par
We find also conditions on the spaces $~E_0,E_1,
F_0 \hbox{ and } F_1~$, so that the following holds
$$
[E_0\wcheck\otimes F_0 , E_1\wcheck\otimes F_1 ]_\th =
[E_0 ,E_1]_\th\wcheck\otimes [F_0 ,F_1]_\th,
\qquad 0 < \th < 1.$$ 
\par
Some applications of these results are also considered.
\bg
\noindent{\bf 1. Introduction, notation and background }
\bg
All Banach spaces considered in this paper are
complex.  
By an $~n$-dimensional Banach space, we mean $\comp^n $ 
equipped with a norm.
\par
If $X$ and $Y$ are Banach spaces, then ${\cal L}(X,Y) $,  $ X\wcheck\otimes Y $ and $ X\hot Y $ denote, respectively, the Banach space of bounded operators 
from $ X $ into $ Y $, The closure of $ X \otimes Y $ 
in  ${\cal L}(X^\ast,Y) $ equipped with the induced
norm, and the completion of $ X \otimes Y $ with 
respect to the projective tensor norm defined by :
$$
\forall u \in X\otimes Y, \qquad\N{u}_\wedge 
= \inf \left\{\sum_{k=1}^m \N{x_k}_X
\N{y_k}_Y : u =  \sum_{k=1}^m x_k\otimes
y_k \right\} $$ 
$ X\wcheck\otimes Y $ and $
X\hot Y $ are called respectively 
the injective  and projective  tensor  
product  of  $X$ and $Y$.  In the
case  when  $X$ and $Y$  are  both  
finite-dimensional, we have  $$
(X\hot Y)^\ast = {\cal L}(X,Y^\ast) = 
 X^\ast\wcheck\otimes Y^\ast. $$  
Using  the preceding  duality,  we see that the results 
announced  in the abstract are similar,  in the  
finite-dimensional  context.
\par
Let $E_0,E_1,F_0 \hbox{ and } F_1$ be
finite-dimensional Banach  spaces.  The usual
interpolation  theorem  asserts that 
$$
\N{.}_{{\cal L}(E_\th^{},F_\th^\ast)} 
\leq
\N{.}_{[{\cal
L}(E_0^{},F_0^\ast),{\cal
L}(E_1^{},F_1^\ast)]_\th} ,$$
 where $ X_\th
$ denotes  $ [X_0,X_1]_\th $ for  $ 0 < \th
< 1 $.   \par
 The question  we  are interested in,  is
the  following: Under what  conditions 
on  the spaces  can one find  a constant
$c$,  independent of the  dimension of the 
 considered spaces,  so that 
$$
\N{.}_{[{\cal
L}(E_0^{},F_0^\ast),{\cal L}(E_1^{},F_1^\ast)]_\th} 
\leq c \N{.}_{{\cal L}(E_\th^{},F_\th^\ast)}.  $$
\par
We will see that  the  constant $c$  can  be  majorized 
 using the  type 2 constants of the spaces  $E_0,E_1,F_0
$ and $ F_1$  ( or the  2-convexity
constants in  the  Banach lattice  case).
\par 
 We  first recall some definitions
and notation. 
${\cal E}_{n,m} $ denotes  the space  of complex 
$ n \times m $ matrices. $ A^\ast $ and $ {}^tA $ are, 
respectively, the  adjoint and the transposed matrix 
of $ A $. A matrix $ A \in {\cal E}_{n,m} $ will be 
identified with a linear operator from $ \comp^m $ 
into $ \comp^n $ using the canonical bases.  
\par 
Let $ \delta $ be a norm on ${\cal E}_{n,m}
$,   the dual  norm is defined by 
$$
\delta^\ast(A) = \sup \left\{ \abs{ \hbox{ tr
}({}^tB.A)} 
 : B \in {\cal E}_{n,m}\quad \hbox{ and }\quad \delta(B)
\leq 1 \right \}. $$
\par
Let ${\cal D } $ , $\adh{\cal D } $ and $
\partial{\cal D } $ denote, respectively the open unit 
disc in $ \comp $, its closure  and its boundary. For 
 $ z \in {\cal D } $  and $ t \in \partial{\cal D } $
we denote by $ P^z(t) $  the Poisson kernel: 
$$ P^z(t) = {{1-\abs{z}^2}\over {\abs{t-z}^2}} $$ 
and let $ dm $ be the Haar  measur  on
$\tore \equiv \partial{\cal D } $.
\par
We use the notation $ A({\cal D }, {\cal E}_{n,m}) $ 
(resp. $ H^{\infty}({\cal D }, {\cal E}_{n,m}) $ )  
to  denote the  set of  analytic functions on ${\cal D
}$  valued in ${\cal E}_{n,m} $  which are
continuous    on  $ \adh{\cal D } $, (resp. bounded 
on  $ {\cal D } $ ).
\par
We say that  an operator $ u : X\rightarrow Y $ 
is $p$-summing for some $ p \geq 1 $ if there is a 
constant  $ c $ such that for all finite  sequences  
$ x_1,...,x_n $ in $ X $  we have 
$$
\left (\sum_{k=1}^n \N{u(x_k)}_Y^p\right )^{1/p} \leq
 c \sup \left\{ \left (\sum_{k=1}^n \abs{\xi(x_k)}^p 
\right )^{1/p} :\quad \xi \in X^\ast ,\quad\N{\xi}_{X^\ast} 
\leq 1 \right \}. $$
We  denote by $\pi_p(u : X\rightarrow Y) $ the smallest 
constant $ c $ satisfying this property, and  by 
$ \Pi_p(X,Y) $ the space of all $p$-summing operators 
 from $ X $ into $ Y $. This space equipped  with 
the $p$-summing norm $ \pi_p $ is a Banach space.
\par
We denote  by $\Gamma_2(X,Y) $ the  set of all
operators which factor through a Hilbert space. For an 
operator $ u \in \Gamma_2(X,Y) $, we define the norm
$$
\gamma_2(u : X\rightarrow Y) = \inf \left\{ 
\N{A :X\rightarrow H}\N{B :H\rightarrow Y} \right\} $$
where the infimum is taken  over all  factorizations 
of $ u $ of the form $ u = B.A $ and all Hilbert spaces
$ H $. The space $\Gamma_2(X,Y) $ epuipped with the  
 preceding norm is a Banach space. 
\par
Let $ 1 \leq p \leq +\infty $. We denote by $ \ell_p^n
$   the  space $ \comp^n $ equipped  with the norm  
$ \N{x}_p = (\sum_1^n \abs{x_k}^p )^{1/p}$, and  with 
the habitual change for $ p = +\infty $. 
\par
Let $ \{ g_n\}_{n\geq 1} $ be a sequence of independent  
identically distributed   Gaussian  real valued normal 
random variables  on some probability space $
(\Omega,{\cal A},\pe )$.  
A Banach space $ X $  is  called of type 2 (resp.
Gaussian cotype 2 ), if there   is a constant $ c $ such
that  for 
 all $ x_1,...,x_n $ in  $ X $ we have 
$$
\N{\sum_{k=1}^n g_k x_k }_{L^2(\Omega;X)} \leq c 
\left ( \sum_{k=1}^n \N{x_k}^2 \right )^{1\over 2} $$ 
(resp. $\geq c^{-1} (\sum \N{x_k}^2)^{1/2} $).
\par
We denote by $ \widetilde T_2(X) $ (resp. $
\widetilde C_2(X) )$  the smallest  constant
$ c $  for which this holds. \par
Let $ (e_1,...,e_n) $ be the canonical basis of 
$ \comp^n $. For any  operator $ u : \ell_2^n 
\rightarrow X $ we define 
$$
\ell(u : \ell_2^n \rightarrow X ) = 
\N{ \sum_{k=1}^n g_k u(e_k) }_{L^2(\Omega;X)}. $$
\par
For more details on operator ideals and  the geometry 
 of Banach  spaces, we refer the reader  to 
[{\bf P}], [{\bf Pi1}] and [{\bf LT1}].
\par
Concerning Banach lattices  we refer the reader to 
 [{\bf LT2}].  We only  recall the following definition: 
   A Banach lattice $ X $ is called 2-convex (resp. 
2-concave ), if  there exists a constant $ M $ so that
for every choice of vectors $ x_1,...,x_n $ in $ X $
$$
\N{\left (\sum_{k=1}^n \abs{x_k}^2\right )^{1/2 }} 
\leq M \left ( \sum_{k=1}^n \N{x_k}^2\right )^{1/2}
$$
(resp. $ \geq M^{-1} (\sum \N{x_k}^2)^{1/ 2}$). The  
smallest possible value of $ M $ is denoted $ 
M^{(2)}(X) $ 
(resp. $ M_{(2)}(X) $). 
\par
We assume the reader familiar with the complex
interpolation  methode of Calder$\acute{\hbox{o}}$n, see 
for   instance [{\bf C}] or [{\bf BL}].
\par
We recall some definitions  and  results  concerning  
the  interpolation  of  families of finite-dimensional 
 Banach  spaces.
\par
A family  of norms  $ \{\delta_\lambda\}_{\lambda \in
\Lambda} $ , where  $ \Lambda $ is any topological
space, on $ {\cal E}_{n,m} $ is said to be  
{\it continuous} 
if for every $ A\in {\cal E}_{n,m} $  the function  
$ \lambda \mapsto \delta_\lambda (A) $ is continuous.
\par
A continuous  family of norms $ \{\delta_z\}_{z\in 
\adh {\cal D} } $  on $ {\cal E}_{n,m} $  is said  to
be    {\it subharmonic}  if for every 
$ {\cal E}_{n,m} $-valued  analytic  function $~F~$
defined on some domain $ \Omega \subset {\cal D} $ 
the function $ z\mapsto \hbox{ Log }\delta_z(F(z)) $ 
is subharmonic  on  $ \Omega $. 
\par
If $ \{\delta_t\}_{t\in \partial {\cal D} } $ is  a  
continuous  family  of norms on $ {\cal E}_{n,m} $,  
there  exists  a  unique  continuous  family  of  
norms  $ \{\delta_z\}_{z\in \adh {\cal D} } $  on  
$ {\cal E}_{n,m} $ which  coincides  with  the 
original  one on the  boundary  and  such that both 
$ \{\delta_z\}_{z\in \adh {\cal D} }$ and  
$ \{\delta_z^\ast\}_{z\in \adh {\cal D} }$ are
subharmonic.  This  family ( or the family  of spaces 
 $ \{ ( {\cal E}_{n,m},\delta_z )\}_{z\in \adh{\cal D}
}$) 
 is called  the interpolation  family with  boundary 
 data 
$ \{\delta_t\}_{t\in \partial {\cal D} } $.
\par
If $ \{\alpha_z\}_{z\in \adh {\cal D} }$ is  an
interpolation family of norms  on $ \comp^m $,  
$ \{\beta_z\}_{z\in \adh {\cal D} }$  is  a subharmonic  
family  of  norms on $ \comp^n $ and 
$ T \in  A({\cal D } ,{\cal E}_{n,m}) $  then the  
function 
$$ z \mapsto \hbox{ Log }\N{ T(z) : (\comp^m,\alpha_z) 
\rightarrow (\comp^n,\beta_z)  } $$
is subharmonic.
\par
Concerning  the preceding  two  results  and more
details on  the properties interpolation  families  we 
refer  the  reader  to [{\bf CS}],[{\bf CCRSW1}]  and
[{\bf CCRSW2}]. We also refer the reader  to this
last reference 
 for  the connection between  this interpolation  
construction  and  the  complex  interpolation method  
 of  Calder\'on. 
\par
Let us finally describe the organization of this paper.
  In section 2, we prove some results which will be
usefull for our problem, in particular,  we prove  a  
more  general  version  of  a theorem due  to Gilles  
Pisier on the  factorization  of matrix-valued 
analytic   functions [{\bf Pi2}].  In section 3, we 
prove  the announced  results  in the 
finite-dimensional  case, using the theory of
interpolation  for families of Banach spaces. In
section 4, we consider  the infinite-dimensional case.
Finally in section 5, we give  some corollaries  and  
applications of our results.
\bg
\noindent{\bf 2. The factorization theorem } 
\bg
The  following  definition  was  introduced in an  
infinite-dimensional  setting in [{\bf Pi2}]. It reduces  
to the  following  in  the  matrix case :
\par
 D{\sevenrm EFINITION} 2.1. A norm  $\delta$  on  $ {\cal
E}_{n,m}$  is  \underbar{2-convex},  if  for every  
$ A , B $ and  $C$  in $ {\cal E}_{n,m}$  the
following  holds 
$$
C^\ast C \leq B^\ast B + A^\ast A \qquad\Longrightarrow\qquad
\delta^2(C) \leq \delta^2(B) + \delta^2(A). $$
The  next lemma is  easy, and its proof is staightforward 
so we  will  omit  it.
\bg
\proclaim Lemma 2.2. 
Let  $ \{ (\alpha_t , {\cal E}_{m,m}   )\}_{t \in
\Lambda } $  and $ \{ (\beta_t , {\cal E}_{m,n} )\}_{t
\in \Lambda } $  be  two continuous  families,  on 
$ \Lambda = \partial{\cal D} \hbox{ or } \adh{\cal D} $,
 formed of 2-convex norms  on  the  corresponding 
matrix  spaces. For  $ t \in \Lambda $ and $ T \in 
{\cal E}_{n,m} $, we  define 
$$
\gamma_t(T) = \inf \left \{ \alpha_t(A).\beta_t(B) :~~~~
 T={}^tB.A ~,~A \in {\cal E}_{m,m}~\hbox{ and }~
B \in {\cal E}_{m,n} \right \}.
$$
Then $ \{ \gamma_t \}_{t \in \Lambda } $  is a
continuous   family of norms  on $ {\cal E}_{n,m} $.
\par
In  particular,  the  fact that $ \gamma_t $ is a norm,
 is a consequence  of the hypothesis  of 2-convexity 
made  on the considered  norms.  We can  now  state  and  
prove  the  main  theorem  of this section. 
\par
\bg
\proclaim Theorem 2.3.  
Let  $ \{ (\alpha_z , {\cal E}_{m,m}   )\}_{z \in
\adh{\cal D} } $  and $ \{ (\beta_z , {\cal E}_{m,n} 
)\}_{z \in \adh{\cal D} } $  be  two  subharmonic 
families  of 2-convex norms. Then the family 
$ \{ (\gamma_z , {\cal E}_{n,m} )\}_{z \in
\adh{\cal D} } $,  defined by 
$$
\gamma_z(T) = \inf \left \{ \alpha_z(A).\beta_z(B) :~~~~
 T={}^tB.A ~,~A \in {\cal E}_{m,m}~\hbox{ and }~
B \in {\cal E}_{m,n} \right \} ,
$$
is subharmonic. More precisely, for every  $ F \in A({\cal
D},{\cal E}_{n,m}) $  and every $ \eps > 0 $ there  exist 
 $ A \in H^\infty({\cal D},{\cal E}_{m,m}) $ and 
 $ B \in H^\infty({\cal D},{\cal E}_{m,n}) $ such that   
for every  $ z \in {\cal D} $ we have 
 $$ F(z) = {}^tB(z).A(z) $$ 
and
$$
\alpha_z(A(z)).\beta_z(B(z)) \leq (1+\eps) \exp 
\int_\tore \hbox{ Log }\gamma_t(F(t)) ~P^z(t) \, dm(t). $$
\par
\bg
P{\sevenrm ROOF}: 
It is clear that we only have to prove  the  last 
assertion  of the theorem.
\par
Let $ F \in A({\cal D},{\cal E}_{n,m}) $  and  $ \eps
> 0 $. There exist  two continuous  functions 
 $ V:\partial{\cal D} \rightarrow {\cal E}_{m,m} $ and 
 $ W:\partial{\cal D} \rightarrow {\cal E}_{m,n} $  
such  that for every $ t \in  \partial{\cal D} $ we have 
$$F(t)~=~{}^tW(t).V(t)  , \eqno (2.1)  $$
$$\alpha_t(V(t))~\leq~1 + {\eps\over 2} , \eqno (2.2)  $$
$$\beta_t(W(t))~\leq~\gamma_t(F(t)).\eqno (2.3) $$
\par
By  the Wiener-Masani theorem [{\bf H},Lecture XI],  
there exists  an  outer  function 
$ A \in H^\infty({\cal D},{\cal E}_{m,m}) $ such that  for  
almost  every  $ t \in \partial{\cal D} $
$$
A^\ast (t)A(t)  = V^\ast (t)V(t) + {{\eps^2}\over{4s^2}} I,
$$
where I is the identity  matrix and  $ s = \sup \{
\alpha_t(I)~:~t \in \partial{\cal D} \} $.
\par
Taking into account (2.2)  and the 2-convexity of 
$ \alpha_t $ we obtain that for almost all 
$ t \in \partial{\cal D} $
$$
\alpha_t(A(t)) \leq \left ( \alpha_t^2(V(t)) +
{{\eps^2}\over 4} \right )^{1/2} \leq 1 + \eps. \eqno
(2.4) $$
\par
On the other hand, since $ V^\ast (t)V(t) \leq
A^\ast (t)A(t) ~~dm \hbox{ a.e.} $  there exists  a 
measurable function 
$ S:\partial{\cal D} \rightarrow {\cal E}_{m,m} $  such
that 
$$
S^\ast (t)S(t) \leq I ~~\hbox{ and }~~V(t) = S(t)A(t) 
 ~~dm \hbox{ a.e.} \eqno (2.5) $$
\par
We can now define  $ B(z) = [\,{}^tA(z)]^{-1}.{}^tF(z) $
for  every  $ z \in {\cal D} $. Clearly $ B $ is an
analytic  function  on $ {\cal D} $ which  is also  
bounded  since  by (2.1) , (2.5) and (2.3) we have 
$$
B(t) = [{}^tA(t)]^{-1} {}^tV(t) W(t) = {}^tS(t)W(t) 
 ~~dm \hbox{ a.e.} $$
$$
\beta_t(B(t)) \leq \beta_t(W(t)) \leq \gamma_t(F(t))  
~~~~~~~~~~dm \hbox{ a.e.}  \eqno (2.6) $$
\par
By the subharmonicity  of the norms  and 
(2.4) , (2.6) we obtain that for every  $ z \in {\cal D} $
 $$\eqalign {&F(z)={}^tB(z).A(z) , \cr
&\alpha_z(A(z)) \leq 1+\eps , \cr
&\beta_z(B(z)) \leq \exp 
\int_\tore \hbox{ Log }\gamma_t(F(t))~ P^z(t) \, dm(t).
}$$
 This achieves the proof of the theorem.
\par
{\it Remark.}
Note that if, in  the preceding theorem, we take 
$ \alpha_z = \alpha~,~\beta_z = \beta $ for all $ z \in 
\adh{\cal D} $, we obtain  the  main part of a  
factorization theorem due to G. Pisier [{\bf Pi2}]. 
\par
\proclaim Corollary 2.4. 
Let $ \{(\alpha_z ,{\cal E}_{m,m})\}_{z\in \adh {\cal D}}$ and
 $ \{(\beta_z ,{\cal E}_{m,n})\}_{z\in \adh {\cal D}}$  be two
interpolation families of 2-convex norms.
Then the family $ \{(\gamma_z,{\cal E}_{n,m})\}_{z\in \adh {\cal
D}}$, defined by 
 $$ \gamma_z(T) = \inf \{ \alpha_z(A).\beta_z(B):~~~T={}^tB.A,
A \in {\cal E}_{m,m} \hbox{ and } B \in {\cal E}_{m,n}\},$$
is an interpolation family.
\par
P{\sevenrm ROOF}: Using the preceding theorem we know that the family 
$ \{(\gamma_z,{\cal E}_{n,m})\}_{z\in \adh {\cal D}}$  is
subharmonic.
\par
Let $ \{( N_z,{\cal E}_{n,m})\}_{z\in \adh {\cal D}}$ be a
subharmonic family such that for every $ t \in \partial {\cal D}
$  and for every $ T \in {\cal E}_{n,m} $ we have $ N_t(T) \le
\gamma_t(T) $. Consider $ T \in {\cal E}_{n,m} $ and $ \omega \in 
{\cal D} $, by definition, there exist $ A \in {\cal E}_{m,m} $
and $ B \in {\cal E}_{m,n}$, such that 
$$ T = {}^tB.A \qquad\hbox{ and }\qquad \alpha_\omega(A) =
\beta_\omega(B) = \sqrt{\gamma_\omega(T)}. $$
\par
Since $ \{\alpha_z\}_{z\in \adh {\cal D}}$ and
 $ \{\beta_z\}_{z\in \adh {\cal D}}$  are interpolation
families, there exist analytic matrix-valued functions 
$ F \in H^\infty({\cal D},{\cal E}_{m,m}) $ and 
$ G \in H^\infty({\cal D},{\cal E}_{m,n}) $, such that 
$$ F(\omega) = A,~G(\omega) = B~~\hbox{ and }~~ \alpha_t(F(t)) =
\alpha_\omega(A) ,~\beta_t(G(t)) = \beta_\omega(B)~~dm~ a.e. $$
\par
By the subharmonicity of the family
$ \{ N_z\}_{z\in \adh {\cal D}}$, we obtain  
$$\eqalign { N_\omega(T) &= N_\omega({}^tG(\omega).F(\omega)) \cr
&\le \int_{\partial {\cal D}} N_t({}^tG(t).F(t))~P^z(t) \,
dm(t) \cr &\le \int_{\partial {\cal D}}
\gamma_t({}^tG(t).F(t))~P^z(t) \, dm(t) \cr
&\le \int_{\partial {\cal D}} \alpha_t(F(t))
\beta_t(G(t))~P^z(t) \, dm(t) \cr
&= \gamma_\omega(T). }$$
\par
We have proved that
$$ \forall~\omega \in {\cal D} , \quad \forall~T \in {\cal
E}_{n,m} , \qquad  \gamma_\omega(T) = \sup \{~ N_\omega(T)~\} $$
where the supremum is taken over the set of  all subharmonic
families $ \{( N_z,{\cal E}_{n,m})\}_{z\in \adh {\cal
D}}$  of norms satisfying the boundary condition: $$
\forall~t \in \partial{\cal D}, \forall~T \in {\cal
E}_{n,m},\qquad N_t(T)~\le~\gamma_t(T). $$ This proves that  
$ \{(\gamma_z,{\cal E}_{n,m})\}_{z\in \adh {\cal D}}$ is an
interpolation family.(\cf [{\bf CS}]).
\par
The next result expresses the stability of the 
property of 2-convexity under the interpolation
construction.
\par
\proclaim Theorem 2.5. 
Let $ \{ (\delta_t ,{\cal E}_{n,m})\}_{t \in
\partial{\cal D} } $ be a continuous family of 2-convex  
 norms. Then the interpolation family $ \{(\alpha_z ,
{\cal E}_{n,m}   )\}_{z \in \adh{\cal D} } $ obtained
from  the preceding one consists of 2-convex norms.
\par
P{\sevenrm ROOF}: 
For $ z \in {\cal D} $ and $ A \in {\cal E}_{n,m}$ define
$$\eqalign {
\Delta_z(A) = \inf \{ \lambda>0 : \exists F \in 
H^\infty({\cal D},{\cal E}_{n,m})~ &\hbox{ with } A^\ast 
A \leq  \lambda^2  F^\ast (z)F(z) \cr
&\hbox{ and }
\delta_t(F(t))\leq 1 ~dm \hbox{ a.e} \}. }$$ 
\proclaim Claim. For every $ z $ in $ {\cal D} $ we have $
\Delta_z = \delta_z $.
\par
Indeed, if $ \delta_z(A) = 1 $ then by Theorem I. from 
[{\bf CCRSW1}] we can find  an element $ F \in H^\infty({\cal
D},{\cal E}_{n,m}) $ such that  $ F(z) = A $  
and  $ \delta_t(F(t)) \leq 1 ~~dm \hbox{ a.e} $, hence 
$ \Delta_z(A) \leq 1 $.
\par
Conversly, if $ \Delta_z(A) < 1 $ there exists an
analytic function $  F \in H^\infty({\cal
D},{\cal E}_{n,m}) $ such that $ A^\ast A \leq F^\ast 
(z)F(z) $ and $ \delta_t(F(t))\leq 1 ~dm \hbox{ a.e.} $ 
 Thus  there exists a contraction $ T $ ,\ie $ T^\ast T
\leq I $, such  that  $ A = T F(z) $. For $ \omega \in 
{\cal D} $ let us consider $ G(\omega) = T F(\omega) $.  
 Clearly $ G \in H^\infty({\cal D},{\cal E}_{n,m}) $ , 
 $ G(z) = A $ and for almost all $ t \in \partial{\cal
D} $ we have  $ G^\ast (t)G(t) \leq F^\ast (t)F(t) $,  
 so $ \delta_t(G(t)) \leq \delta_t(F(t)) \leq 1 $. 
 This proves the claim.
\par 
 Let us consider $ A, B \hbox{ and } C $ in $ {\cal
E}_{n,m} $ and  $ \eps > 0 $.
\par
Define $ \alpha = \delta_z(A) $ and $ \beta =
\delta_z(B) $. Without loss of generality,  we may
assume  that $ \alpha^2 + \beta^2 = 1 $. Then there 
exist $ F $ and $ G $ in  $ H^\infty({\cal D},{\cal
E}_{n,m}) $ such that  $ F(z) = A , G(z) = B $ and  for  
almost all  $ t \in \partial{\cal D}~~~\delta_t(F(t)) =
\alpha , \delta_t(G(t)) = \beta $.
\par
By the Wiener-Masani theorem [{\bf H},Lecture XI], there 
exists an outer function $ L \in H^\infty({\cal D},{\cal
E}_{n,m}) $ such that 
$$
L^\ast (t)L(t) = F^\ast (t)F(t) + G^\ast (t)G(t) + 
{{\eps^2}\over {s^2}} I ~~~dm \hbox{ a.e.} \eqno (2.7) $$
with $ s = \sup \{ \delta_t(I) ~:~t \in  \partial{\cal
D }\}$.
\par
By (2.7) and using the 2-convexity of $ \delta_t $ we
obtain 
$$
\delta_t(L(t)) \leq ( \alpha^2 + \beta^2 + \eps^2 )^{1/2} 
\leq 1 + \eps  ~~~dm \hbox{ a.e.} \eqno (2.8) $$
\par
Take $ x \in \comp^m $ and $ h \in H^\infty({\cal
D},\comp^m) $ such that  $ h(z) = x $, and consider the 
function $ k \in H^\infty({\cal D},\comp^{2n}) $ defined 
by
$$
k(\omega) = \left[ \matrix{F(\omega) h(\omega)\cr {}\cr
G(\omega) h(\omega)}\right].$$
\par
Since $ k $ is  analytic, the function  $ \omega \mapsto 
\hbox{ Log }\N{k(\omega)}_2 $ 
is subharmonic, (where $ \N{.}_2 $ is the Hilbertian
norm), and we may write: 
$$\left(\N{Ax}_2^2 +\N{Bx}_2^2 \right)^{1/2} = 
\N{k(z)}_2  \leq \exp 
\int_\tore \hbox{ Log }\N{k(t)}_2~ P^z(t) \, dm(t).
\eqno (2.9) $$
By (2.7) ,(2.9) and the hypothesis we obtain
$$
\N{ Cx}_2 \leq \exp 
\int_\tore \hbox{ Log }\N{L(t)h(t)}_2~ P^z(t) \, dm(t).
\eqno (2.10) $$
By the  results of [{\bf CCRSW1}] on the interpolation 
of Hilbertian families of norms, we know that  the 
family  $ \{ y\mapsto \N{ L(\omega)y }_2 \}_{\omega
\in \adh {\cal D} } $  is an interpolation family of 
norms  on $ \comp^m $. Consequently, by taking the infimum 
in (2.10) over all $ h \in H^\infty({\cal D},\comp^m) $
such that  $ h(z) = x $  we obtain
$$
\N{Cx}_2  \leq \N{L(z)x}_2 $$
but $ x \in \comp^m $ is arbitrary, so we have 
$$ C^\ast C \leq L^\ast (z)L(z). \eqno (2.11) $$
Using (2.8) and (2.11) we obtain that $ \Delta_z(C) 
\leq 1+\eps $ and the result follows since $ \eps $ is 
arbitrary.
\par
\bg
\noindent{\bf 3. The finite-dimensional case }
\bg
This section is devoted to the proof of the following
theorem.
\par
\proclaim Theorem 3.1. 
Let $ \{ E_z\}_{z \in \adh{\cal D}} $ and  $ \{ F
_z\}_{z \in \adh{\cal D}} $ be two interpolation families 
of $ m$-dimensional and $ n$-dimensional spaces
respectively. 
If $ T \in {\cal E}_{n,m} $ and $z \in \adh{\cal D} $, 
define $ \NN{T}_z $ as the operator norm : $
\N{T : E_z \rightarrow F_z^\ast} $ and  
$\{ \N{.}_{[z]} \}_{z \in \adh{\cal D}} $ as the
interpolation family of norms on ${\cal E}_{n,m} $
extending $ \{ \NN{.}_t\}_{t \in \partial {\cal D}} $. 
 Then 
$$\forall~ z\in \adh{\cal D}, \forall~T \in {\cal E}_{n,m}
~~~~~\NN{T}_z \leq \N{T}_{[z]} \leq c \NN{T}_z , $$
with the following estimates for c :\hbi
1. $ c \leq \left(\sup~ \{
\widetilde T_2(E_t).\widetilde T_2(F_t)~ :~t \in \partial
{\cal D} \} \right)^2 $. \hbi
2. $ c \leq \displaystyle{2\sqrt{{2\over \pi}}}\left(\sup~ \{
\widetilde T_2(E_t).M^{(2)}(F_t)~ :~t \in \partial {\cal D} \}
\right)^2 $ if  the canonical basis of $ \comp^n$ 
 is a 1-unconditional basis for $ F_t $ for every
$t \in \partial {\cal D} $.\hbi
3. $ c \leq \displaystyle{{32\over\pi}}\left(\sup~
\{M^{(2)}(E_t).M^{(2)}(F_t)~ :~t \in \partial {\cal D} \}
\right)^{5/2} $ if  the canonical basis of $ \comp^m$ 
(resp.$ \comp^n$ ) is a 1-unconditional basis for $ E_t $
(resp.$ F_t $ ) for every $t \in \partial {\cal D} $.
\par
\bg
The proof of this theorem uses several lemmas which will
be listed and proved in what follows.
\proclaim Lemma 3.2. 
Let $ \{ E_z\}_{z \in \adh{\cal D}} $   be an
interpolation family  of $ m$-dimensional  spaces.  If $ T
\in {\cal E}_{n,m} $ and $z \in \adh{\cal D} $,  define 
$\pi_z^E(T) = \pi_2(T: E_z^\ast \rightarrow \ell_2^n) $, 
 and let $ \{ \pi_{[z]}^E \}_{z \in \adh{\cal D}} $  be 
the interpolation family of norms on $ {\cal E}_{n,m} $
extending 
$ \{ \pi_t^E \}_{t \in \partial{\cal D}} $. Then
$$\forall~ z \in \adh{\cal D} , \forall~  T
\in {\cal E}_{n,m} ~~~~\pi_z^E(T) \leq W^E(z)
\pi_{[z]}^E(T) \eqno (\ast) $$
with $ W^E(z) = \exp \int_\tore \hbox{ Log } \wt
T_2(E_t) ~P^z(t) \, dm(t) $.
\par
\bg
P{\sevenrm ROOF} of lemma 3.2. 
 Note that if $ X $ is $ m$-dimensional then for every $ T
\in {\cal E}_{n,m} $ one has 
$$
\pi_2(T: X^\ast \rightarrow \ell_2^n)  \leq~ 
\ell ({}^tT:\ell_2^n\rightarrow X) \leq~ 
\widetilde T_2(X). \pi_2(T: X^\ast \rightarrow \ell_2^n).
\eqno (3.1) $$ 
\par
The first inequality follows trivially from the
definition,  and the second one is well-known and
classical, see for instance [{\bf P},21.3.5].
\par
Define $ \ell_z(T) = \ell ({}^tT:\ell_2^n\rightarrow E_z)
$  for $ z \in \adh{\cal D} $   and  $ T\in {\cal E}_{n,m}
$.  
Clearly $ \{ \ell_z\}_{z \in \adh{\cal D}} $ is a
subharmonic family of norms on $ {\cal E}_{n,m} $.  
Using the result of [{\bf CS}] recalled in section 1. we 
obtain that the function 
$$
z \longmapsto \hbox{ Log}\N{ I:({\cal
E}_{n,m},\pi_{[z]}^E) \rightarrow ({\cal
E}_{n,m},\ell_z)} $$ is subharmonic, so by (3.1) we may
write 
$$
\forall~ z \in \adh{\cal D}, ~~\N{ I:({\cal
E}_{n,m},\pi_{[z]}^E) \rightarrow ({\cal
E}_{n,m},\ell_z)} \leq W^E(z).$$
A second  use of (3.1) yields the result.
\par
\proclaim Lemma 3.3.
Let $ \{ E_z\}_{z \in \adh{\cal D}} $   be an
interpolation family  of $ m$-dimensional  spaces, and 
assume that the canonical basis of $ \comp^m $ is a 
1-unconditional  basis for $ E_t $ for every $ t \in 
\partial{\cal D}$. If  $ \{ \pi_z^E \}_{z \in \adh{\cal
D}} $ and $ \{ \pi_{[z]}^E \}_{z \in \adh{\cal D}} $  are 
defined as in lemma 3.2 then
$$\forall~ z \in \adh{\cal D} , \forall~   T
\in {\cal E}_{n,m} ~~~~\pi_z^E(T) \leq W^E(z)
\pi_{[z]}^E(T) \eqno (\ast) $$
with $ W^E(z) = \displaystyle{{2\over{\sqrt{\pi}}}}\exp
\int_\tore\hbox{ Log } M^{(2)}(E_t) ~P^z(t) \, dm(t) $.
\par
\bg
P{\sevenrm ROOF} of lemma 3.3. 
It is well-known ( see [{\bf Pi1},Chapter 8]) that if 
$ X $ is an $m$-dimensional complex Banach lattice, then 
 for every $ T \in {\cal E}_{n,m} $ one has 
$$
{{\sqrt{\pi}}\over 2} \pi_2(T: X^\ast \rightarrow
\ell_2^n) \leq~\N{\left(\sum_{k=1}^n\abs{{}^tT(e_k)}^2
\right)^{1/2}}_X \leq~M^{(2)}(X)~ \pi_2(T: X^\ast
\rightarrow \ell_2^n) \eqno (3.2) $$
\par
 For $ T \in {\cal E}_{n,m}$ and 
$z \in \adh{\cal D} $, we define 
$$\ell_z(T) = \N{\left(\sum_{k=1}^n\abs{{}^tT(e_k)}^2
\right)^{1/2}}_{E_z} = \N{\{{}^tT(e_k)\}_{1\leq k\leq
n}}_{E_z(\ell_2^n)} $$
By the results of [{\bf He}] on the interpolation of
families of Banach lattices,  $ E_z $ is an
$m$-dimensional complex Banach lattice and  
$\{ E_z(\ell_2^n) \}_{z \in \adh{\cal D}} $  
is an interpolation family. Consequently, 
$\{ \ell_z\}_{z \in \adh{\cal D}}$ is subharmonic, and
the proof can be completed as in lemma 3.2 using (3.2) 
instead of (3.1).
\par
{\it Remark.}
Note that using the same notation as in lemma 3.2 and
 lemma 3.3,  we also have  
$$
\pi_{[z]}^E(T) \leq \pi_z^E(T) $$
for all $ z \in \adh{\cal D} $ and all  $ T \in 
{\cal E}_{n,m}$. This will not be used in the sequel.
\par
\bg
\proclaim Lemma 3.4. 
Let $ \{ E_z\}_{z \in \adh{\cal D}} $ and  $ \{ F
_z\}_{z \in \adh{\cal D}} $ be two interpolation families 
of $ m$-dimensional and $ n$-dimensional spaces
respectively. Assume that both families satisfy the 
conclusion $ (\ast) $ of lemma 3.2 or lemma 3.3. 
For every $z \in \adh{\cal D} $ and every $ T \in {\cal 
E}_{n,m}$  define 
$ \gamma_z(T) = \gamma_2( T : E_z \rightarrow F_z^\ast )
$ and  let $ \{ \gamma_{[z]} \}_{z \in \adh{\cal D}} $
  be the interpolation family of norms on $ {\cal
E}_{n,m} $ extending $ \{ \gamma_t\}_{t \in \partial{\cal
D}} $.  Then
$$
\forall~ z \in \adh{\cal D} ,~~\forall~ T \in {\cal E}_{n,m}
~~\gamma_z(T) \leq \gamma_{[z]}(T) \leq K(z) \gamma_z(T)
$$
with  $ K(z) = W^E(z) W^F(z) $.
\par
\bg
P{\sevenrm ROOF} of lemma 3.4. We keep the notation of lemma 3.2 and
  lemma 3.3. Using a well-known fact, (see [{\bf
Pi1},Chapter 3]  for $ T \in {\cal E}_{n,m}$ and 
$z \in \adh{\cal D} $ we may write 
$$\eqalign { \gamma_z^\ast (T) &= \gamma_2^\ast ( T : 
E_z^\ast \rightarrow F_z ) \cr
&= \inf \left \{ \pi_z^E(A).\pi_z^F(B)~:~ T={}^tB.A,
~A \in {\cal E}_{m,m} ,~B \in {\cal E}_{m,n} \right\}.}$$
By $ (\ast) $ we obtain 
$$
\forall~ z \in \adh{\cal D} ,~~\forall~ T \in {\cal E}_{n,m}
~~\gamma_z^\ast(T) \leq K(z)
M_z(T)   \eqno (3.3) $$ 
with
$$
M_z(T)= \inf \left \{ \pi_{[z]}^E(A).\pi_{[z]}^F(B)~:~
T={}^tB.A, ~A \in {\cal E}_{m,m} ,~B \in {\cal E}_{m,n}
\right\}.$$
\par
Using theorem 2.5 we see easily that 
$ \{ \pi_{[z]}^E \}_{z \in \adh{\cal D}} $ and 
$ \{ \pi_{[z]}^F \}_{z \in \adh{\cal D}} $ are
subharmonic families of 2-convex norms and by theorem 2.3 
we see that $ \{ M_z \}_{z \in \adh{\cal D} } $ is a
subharmonic family of norms on $  {\cal E}_{n,m} $; 
moreover, for every $ t \in \partial {\cal D } $  we have
$ M_t = \gamma_t^\ast $. 
 Interpolation and duality theorems imply that 
$$
\forall~ z \in \adh{\cal D}, \forall~ T \in {\cal E}_{n,m}
~~~~~M_z(T) \leq \left( \gamma_{[z]}\right)^\ast(T). \eqno
  (3.4) $$
From (3.3), (3.4) and by  duality we obtain
$$
\forall~ z \in \adh{\cal D}, \forall~ T \in {\cal E}_{n,m} 
~~~~~\gamma_{[z]}(T) \leq K(z) \gamma_z(T) \eqno (3.5) $$
On the other hand, we have also 
$$
\gamma_z(T) = 
\inf \left \{ \N{A: E_z\rightarrow \ell_2^m}.\N{B:
F_z\rightarrow \ell_2^m}~:~ T={}^tB.A, ~A \in {\cal
E}_{m,m} ,~B \in {\cal E}_{m,n} \right\}
$$
and the families $ \{ \N{.: E_z\rightarrow
\ell_2^m} \}_{z \in \adh{\cal D} }$ and $ \{ \N{.:
F_z\rightarrow \ell_2^m} \}_{z \in \adh{\cal D} }$  
  are clearly subharmonic families  of 2-convex norms.  
 Hence by theorem 2.3  we obtain that $ \{ \gamma_z\}_{z 
\in \adh{\cal D}} $  is a subharmonic family, and the
interpolation theorem implies :
$$
\forall~ z \in \adh{\cal D}, \forall~ T \in {\cal E}_{n,m}
 ~~~~\gamma_z(T) \leq \gamma_{[z]}(T)  \eqno (3.6) $$
This completes the proof of the lemma.
\par
\bg
P{\sevenrm ROOF} of theorem 3.1. 
We assume that we are in one of the cases considered in
the theorem, and we keep the notation of the preceding
lemmas.  
By the interpolation theorem  $ \{
\NN{.}_z\}_{  z \in \adh{\cal D}} $ is a subharmonic
family of norms, hence
   
$$
\forall~ z \in \adh{\cal D}, \forall~ T \in {\cal E}_{n,m}
~~~~ \NN{T}_z \leq \N{T}_{[z]}. $$
On the other hand, for every $ t \in \partial{\cal D} $ 
we have clearly 
$$
\N{ I : ({\cal E}_{n,m},~\gamma_t) \rightarrow 
({\cal E}_{n,m},~ \NN{.}_t) } \leq 1.$$
So by interpolation, we conclude that
$$
\forall~ z \in  \adh{\cal D}~~~~
\N{ I : ({\cal E}_{n,m},~\gamma_{[z]}) \rightarrow 
({\cal E}_{n,m},~ \N{.}_{[z]}) } \leq 1.$$
Since lemma 3.4 applies in the cases considered in
theorem 3.1, we have
$$
\forall~ z \in \adh{\cal D}, \forall~ T \in {\cal E}_{n,m}
~~~~\N{T}_{[z]} \leq  \gamma_{[z]}(T) \leq
K(z).\gamma_z(T) \eqno (3.7) $$
\par
Using Kwapie\'n's theorem [{\bf Kw}] or 
[{\bf Pi1},Chapter 3],  in the first case of 
theorem 3.1 we have 
$$
\gamma_z(T) \leq \widetilde T_2(E_z) \widetilde T_2(F_z)
\NN{T}_z  \eqno (3.8) $$
since the Gaussian cotype 2 constant of $ F_z^\ast $ is
smaller than $ \widetilde T_2(F_z) $.
\par
By the same theorem, in the second
case   of theorem 3.1 we have 
$$
\gamma_z(T) \leq \sqrt{2}~ \widetilde T_2(E_z) M^{(2)}(F_z)
\NN{T}_z  \eqno (3.9) $$
because the Gaussian cotype 2 constant of $ F_z^\ast $ is 
smaller than $ \sqrt{2} M^{(2)}(F_z) $.
\par
 Here we used the simple fact that when we
interpolate a family of 2-convex Banach lattices 
$ \{ F_t \}_{t \in \partial{\cal D}} $, we obtain a 
family of 2-convex Banach lattices $ \{ F_z \}_{z 
\in \adh{\cal D}} $. Moreover, the function 
$ z \mapsto \hbox{ Log }M^{(2)}(F_z) $ is subharmonic. 
See [{\bf He }] for more details.
\par
Finally, using theorem 4.1 from [{\bf Pi1}], instead of 
 Kwapie\'n's theorem, we have in the third
case :
$$
\gamma_z(T) \leq \left (4 M^{(2)} (E_z)
M^{(2)}(F_z)\right)^{3/2} \NN{T}_z  \eqno (3.10) $$ 
\par
If $ \{ X_z \}_{z \in \adh{\cal D}} $ is an interpolation 
family of finite-dimensional spaces (resp. Banach
lattices), then the function  $ z \mapsto \hbox{ 
Log }\widetilde T_2(X_z) $ (resp. $ z \mapsto \hbox{
 Log }M^{(2)}(X_z) $ ) is easily seen to be subharmonic. 
 Using this fact, the estimates (3.8), (3.9), (3.10) 
and the estimates of $ K(z) $ obtained in lemma 3.2 and 
lemma 3.3,  we obtain from (3.7) that 
$$
\forall~ z \in \adh{\cal D}, \forall~ T \in {\cal E}_{n,m}
~~~~   \N{T}_{[z]} \leq c \NN{T}_z. $$
where $ c $ is majorized exactly as in the theorem. This
achieves the proof of theorem 3.1.
\par
\bg
\proclaim Corollary 3.5.
Let $ E_0, E_1, F_0 $ and $ F_1 $ be finite-dimensional
Banach spaces with $ \dim E_0 = \dim E_1 = m $, 
$ \dim F_0 = \dim F_1 = n $, and $ \th \in ]0,1[ $. 
Then for all $ T \in {\cal E}_{n,m}$
$$
\N{T}_{{\cal L}(E_\th^{},F_\th^\ast) } \leq 
\N{T}_{[{\cal L}(E_0^{},F_0^\ast),{\cal
L}(E_1^{},F_1^\ast)]_\th} \leq c 
\N{T}_{{\cal L}(E_\th^{},F_\th^\ast) }, $$
with $ E_\th = [E_0,E_1]_\th $ , $ F_\th = [F_0,F_1]_\th $
  and $ c $ is majorized as follows.\hbi
1. $ c \leq \left(\max \{ \widetilde T_2(E_i).\widetilde
T_2(F_i)~ :~i=0,1 \} \right)^2 $.\hbi
2. $ c \leq\displaystyle{ 2\sqrt{{2\over \pi}}}\left(\max \{
\widetilde T_2(E_i).M^{(2)}(F_i)~ :~i = 0,1 \}
\right)^2 $ if  the canonical basis of $ \comp^n$ 
 is a 1-unconditional basis for $ F_i $ for $ i = 0,1
$.\hbi
3. $ c \leq \displaystyle{{32\over\pi}}\left(\max
\{M^{(2)}(E_i).M^{(2)}(F_i)~ :~i = 0,1 \}
\right)^{5/2} $ if  the canonical basis of $ \comp^m$ 
(resp.$ \comp^n$ ) is a 1-unconditional basis for $ E_i $
(resp.$ F_i $ ) for $ i = 0,1 $.
\par
\bg
P{\sevenrm ROOF}: 
This corollary is a simple consequence of theorem 3.1. 
We only need to use a simple approximation argument, and 
theorem 5.1 from [{\bf CCRSW2}] which makes the link
between interpolation families and the complex
interpolation method of Calder\'on.
\par
\bg
\proclaim Corollary 3.6. 
Let $ E_0, E_1, F_0 $ and $ F_1 $ be finite-dimensional
Banach spaces with $ \dim E_0 = \dim E_1 = m $, 
$ \dim F_0 = \dim F_1 = n $, and $ \th \in ]0,1[ $. 
Then for all $ T \in {\cal E}_{n,m}$
$$
C^{-1}(F) ~\pi_2(T : F_\th^\ast \rightarrow E_\th) 
\leq
\N{T}_{[\Pi_2(F_0^\ast,E_0^{})
, \Pi_2(F_1^\ast,E_1^{})]_\th}  \leq C(E^\ast)~ \pi_2(T :
F_\th^\ast \rightarrow E_\th) $$
 with  $ E_\th =
[E_0,E_1]_\th $ , $ F_\th = [F_0,F_1]_\th $;  if $ X $
denotes either $ E^\ast $ or $ F $, $ C(X) $ is majorized as
follows: \hbi
1. $ C(X) \leq \left(\max~ ( \widetilde T_2(X_0),
\widetilde T_2(X_1) )\right)^2 $.\hbi
2. $ C(X) \leq\displaystyle{ 2\sqrt{{2\over \pi}}}\left(\max~
(M^{(2)}(X_0) , M^{(2)}(X_1) ) \right)^2 $ if 
the canonical basis of $ \comp^p~$  $( p = \dim X_i, i =
0,1 ) $ 
 is a 1-unconditional basis for $ X_i $ for $ i = 0,1 $.
\par
\bg
P{\sevenrm ROOF}: 
We reproduce the method of Pisier in [{\bf Pi3}].
Consider 
 the bilinear operator 
$$
\Phi_i : {\cal L}( \ell_2^r , F_i^\ast ) \times 
 \Pi ( F_i^\ast , E_i ) \longrightarrow \ell_2^r(E_i) , $$
where $ i = 0,1 $, defined by 
$ \Phi_i(T , U) = ( U\rond T ( e_k))_{k \leq r} $ with 
 $ r $ any positive integer and $ ( e_k)_{k \leq r} $ 
is the canonical basis of $ \ell_2^r $.
\par
For $ i =0,1 $ we have by definition $ \N{ \Phi_i } \leq
1 $, hence by interpolation we obtain , for all 
$ T \in {\cal E}_{m,r}$ and all  $ U \in {\cal E}_{n,m}$
that 
$$
\left ( \sum_{k=1}^r \N{ U\rond T (e_k)}_{E_\th}^2
\right)^{1/2} \leq 
\N{T}_{[{\cal L}(\ell_2^r,F_0^\ast),{\cal
L}(\ell_2^r,F_1^\ast)]_\th}.
\N{U}_{[\Pi_2(F_0^\ast,E_0^{}),\Pi_2(F_1^\ast,E_1^{})]_\th}. $$
Using corollary 3.5 we obtain
$$ \forall ~r \in \nat ,~~\forall~T \in {\cal E}_{m,r} ,~~
\forall~U \in {\cal E}_{n,m} $$
$$
\left ( \sum_{k=1}^r \N{ U\rond T (e_k)}_{E_\th}^2
\right)^{1/2} \leq C(F)
\N{T}_{{\cal L}(\ell_2^r,F_\th^\ast)}.
\N{U}_{[\Pi_2(F_0^\ast,E_0^{}),\Pi_2(F_1^\ast,E_1^{})]_\th}.
$$ This clearly implies the first inequlity. The second
one  follows by duality.
\par
\bg
{\it Remark. }
 It is shown in the proof of lemma 3.4 that, if $
\{E_z\}_{z \in \adh{\cal D}} $ is an interpolation
family,  then $ \{ \gamma_2(. : E_z \rightarrow E_z) \}_
{z \in \adh{\cal D}} $ is a subharmonic family of norms
on $ {\cal E}_{n,n} $ where $ n = \dim E_z $. Taking the
identity operator and noting that $ 
\gamma_2(I : E_z \rightarrow E_z) = d(E_z , \ell_2^n) $ ; 
the Banach-Mazur distance between $ E_z $ and $ \ell_2^n $
 we infer that the function $ z \mapsto \hbox{ Log }d(
E_z , \ell_2^n) $ is subharmonic. In particular, if 
$ E_0 $ and $ E_1 $ are two $n$-dimensional Banach spaces,
the function  $ \th \mapsto d([E_0,E_1]_\th , \ell_2^n) $
is  Log-convex on the interval $ [0 , 1] $.
\par
\line{}\par
\bg
\noindent{\bf 4. The infinite-dimensional case }
\bg
Let $ X_0 , X_1 $ be two Banach  spaces. The pair $ (X_0
, X_1) $ is  called an interpolation pair if both spaces 
$ X_0 $ and $ X_1 $ are continuously embedded in a
Banach space $ {\cal U } $ so that we can define their
sum and  intersection and we can construct the
interpolation  space $ X_\th = [ X_0 , X_1 ]_\th $ by
 the complex interpolation method of
Calder\'on.
\par
In what follows, all interpolation pairs are assumed to
have dense intersections,\ie $ X_0 \cap X_1 $ is a dense
subspace  of $ X_k,  k = 0 ,1 $.
\par
We will say that $ (X_0 , X_1) $ is an interpolation
pair of Banach lattices if $ X_0 $ and $ X_1 $ are
Banach  lattices which can be embedded as
lattice-ideals  in some $ L^0(\Omega, {\cal A}, \mu ) $, 
with $ (\Omega, {\cal A}, \mu ) $ a $ \sigma$-finite 
measure space.  We assume, moreover, that the set of
simple functions 
$$
S_X = \left \{ \sum_{k=1}^n \alpha_k \wa_{\Omega_k} ~:~
n \in \nat ,~\alpha_k \in \comp ,~\Omega_k \in {\cal A}~ 
\hbox{ and }~ \mu(\Omega_k) < +\infty \right \}
$$
is dense in both $ X_0 $ and $ X_1 $.
\par
If $ (X_0 , X_1) $ is an interpolation pair, and $ G $
is a finite-dimensional subspace of $ X_0 \cap X_1 $, 
then $ G_k $ denotes $ G $ considered as a
subspace of $ X_k$ with $ k= 0,1 $ whereas $ G_\th $
denotes $ [G_0,G_1]_\th $ and not $ G $ considered as a
subspace $ X_\th $.
\par
The following lemma will be the main tool to obtain
infinite-dimensional results from the finite-dimensional
ones, obtained in the preceding section.
\par
This lemma is not new. The author thanks Professor 
A.V. Bukhvalov for indicating to him the book [{\bf KP}]
as a reference for more general results and for telling
him that such results also appear  in Russian literature
in the works of Aizenshtein and Brudnyi,(see [{\bf
Ai}] and [{\bf AB}]). However, we will  include its proof
for the convenience of the reader.
\par
\proclaim Lemma 4.1.
 Let $ (X_0 , X_1) $ be an
interpolation, $ X $ a subspace of $ X_0 \cap X_1 $
which  is dense  both in $ X_0 $ and $ X_1 $, and $ \th 
\in ]0 , 1[ $. Then for every $ \eps \in ]0 , 1[ $  and 
every finite-dimensional subspace $ G $ of $ X $, there
exists a finite-dimensional  subspace $ \widetilde G $ of
$ X $ containing $ G $ and such that 
$$ 
\forall ~x \in G ~~~~~~(1-\eps)~\N{ x }_{\widetilde
G_\th}  \leq~\N{ x }_{X_\th} \leq~\N{ x
}_{\widetilde G_\th}. $$ Moreover, if $(X_0 , X_1) $ is
an interpolation pair of
 Banach lattices embedded in $ L^0(\Omega, {\cal A}, \mu )
$,  and if $ G $ is spanned by finite sequence of
characteristic functions of measurable sets, then $ G $ 
can also be taken of this form.
\par
\bg
P{\sevenrm ROOF}: 
Note that the second inequality is obvious. Let us prove
the first one.
\par
Let $ 0 < \delta < {\eps\over 2} $ and $ {\cal N} = \{
x_1, ...,x_m \} $ be $\delta$-net  in the unit sphere of 
$ (G, \N{.}_{X_\th} ) $. ( \ie for every $ k \leq m $, we
have $ \N{x_k}_{X_\th} = 1 $, and for every $ y \in G $
with $  \N{y}_{X_\th} = 1 $  there exists  $ j \leq m $ 
such that $ \N{ y - x_j }_{X_\th} \leq ~\delta $ ).
Clearly such a net exists since $ G $ is
finite-dimensional.
\par
As $ \N{x_k}_{X_\th} = 1$, for every $ k \in \{
1,...,m\} $ there exists $$
F_k(z) = \sum_{r=1}^{m_k} ~\Psi_{r,k} (z) x_{r,k} $$
with $ x_{r,k} \in X ,~~\Psi_{r,k} \in A(S) $ ( where $ S
$ denotes the stripe $ \{ z~:~0 <\hbox{Re} z < 1 \}$ ) 
such that $ F_k(\th) = x_k $ and for every $ t \in \reel
$ and every $ j = 0,1 $
$$
\N{ F_k( j + it )}_{X_j} \leq~1 + \delta .
$$
\par
Let  $ \widetilde G =~\hbox{ span } \{ x_{r,k}~:~
1 \leq~k \leq~m, 1 \leq~r\leq~m_k \} $. This is clearly
a finite-dimensional subspace of $ X $ containing $ G $.
\par
Consider $ x \in G $ with $ \N{x}_{X_\th} = 1 $, then
one may write $ x = y_0 + \sum_{k=1}^\infty ~\lambda_k
y_k $ with $ 0 \leq~\lambda_k \leq~\delta^k $ 
and 
$ \forall~k, ~~~~y_k \in {\cal N} $. If $ y_k = x_j $ 
we put $ H_k = F_j $ and define $  F = H_0 +
\sum_{k=1}^\infty ~\lambda_k H_k $. We see immediately 
that $ F $ is an analytic $ \widetilde G$-valued  function
on $ S $, continuous on $ \adh S $, satisfying
$ F(\th) = x $ and
$$
\forall~t \in \reel ~,\forall~j= 0,1 ~~~~~
\N{ F(j + it) }_{X_j} \leq~{{1 +
\delta}\over{1 - \delta}}. $$
This implies that 
$$ \N{x}_{\widetilde G_\th} \leq~{{1 +
\delta}\over{1-\delta}} \leq ~{1\over{1 - \eps}}.$$
\par
The assertion concerning Banach lattices is clear since,
under our assumptions, the space $ X = S_X $ of simple
functions is dense  both in $ X_0 $ and $ X_1 $.
\par
\bg
\proclaim Theorem 4.2. 
Let $ ( E_0 , E_1 ) $ and $ ( F_0 , F_1 ) $ be two
interpolation pairs. Assume that one of the
following properties holds
$$\eqalign {~~~1&.~E_0^\ast, ~E_1^\ast, ~F_0^\ast
\hbox{ and } ~F_1^\ast~\hbox{ are type 2 spaces }.\cr
~~~2&.~E_0^\ast, ~E_1^\ast~\hbox{ are type 2 spaces 
and } ~F_0^{}, ~F_1^{}~\hbox{ are 2-concave Banach
lattices}.\cr
~~~3&.~E_0^{}, ~E_1^{}, ~F_0^{}
\hbox{ and } ~F_1^{}~\hbox{ are 2-concave Banach
lattices}. } $$
Then $ ( E_0 \wcheck\otimes F_0 ,E_1 \wcheck\otimes F_1 )$ 
 is an interpolation pair, and for $ 0 < \th < 1 $ we
have
$$
\left[E_0 \wcheck\otimes F_0 ,E_1 \wcheck\otimes F_1\right
]_\th = 
\left[E_0,E_1\right]_\th \wcheck\otimes \left[F_0,F_1\right
]_\th. $$
\par
\bg  
P{\sevenrm ROOF}: 
Note first that $ E_j\wcheck\otimes F_j ~,~j= 0,1 $, is
continuously embedded in $ (E_0 + E_1) \wcheck\otimes
( F_0 + F_1 )$. This is assertion (5) in paragraph 44.4 
of K\"othe's book [{\bf K\"o}]. So $ ( E_0 \wcheck\otimes
 F_0 ,E_1 \wcheck\otimes F_1 )$  is an interpolation pair.
\par
Let us make the following convention of notation. If 
$( X_0 ,X_1 ) $ is an interpolation pair of Banach
lattices, then $ X $ denotes the set of simple functions 
$S_X $ defined in the begining of this section which is  
dense in both $ X_0 $ and $ X_1 $. If $ X_0$ and $X_1 $
are not assumed to be Banach lattices then $ X $ is
simply  $ X_0 \cap X_1 $. With this convention, we 
easily see that $ E \otimes F $ is a dense subspace of 
$ E_\th \wcheck\otimes F_\th $, for all $ \th \in [0,1]
$.  
 In order to prove the result it is sufficient to show
that  $$
\forall~T \in~ E \otimes F~~~
\N{T}_{\wco{E_\th}{F_\th}} \leq~\N{T}_{
[\wco{E_0}{F_0} , \wco{E_1}{F_1}]_\th^{}} \leq~
c \N{T}_{\wco{E_\th}{F_\th}}. \eqno (4.1) $$
\par
If $ ~E_0^{}, ~E_1^{}, ~F_0^{}
\hbox{ and } ~F_1^{} $ are finite-dimensional then 
$ {\cal L}( E_\th^\ast, F_\th^{}) = E_\th \wcheck\otimes 
F_\th $  and the theorem in this case is exactly
corollary 3.5.
\par
In order to prove the infinite-dimensional case, we will
use lemma 4.1. Full details of the argument will be
given only in the first case (\ie when $ ~E_0^\ast, 
~E_1^\ast, ~F_0^\ast
\hbox{ and } ~F_1^\ast$ are spaces of type 2.
 The other cases are treated similarly using the second 
part of lemma 4.1 when  necessary. Details are easy
and left as an exercise for the reader.
\par
We assume that we are under the hypothesis of the first 
case. 
\par
Let $ T = \sum_{k=1}^m~x_k\otimes y_k \in E \otimes F $ 
and $ \eps > 0 $. Define $ X = \hbox{ span }~\{ x_k
~:~k\leq m \} $ and $ Y = \hbox{ span }~\{ y_k ~:~k\leq
m \} $.  By lemma 4.1 we can find $ G $ and $ H $
finite-dimensional subspaces of $ E $ and $ F $
respectively, such that $ X \subset G ~,~Y \subset H $ 
and
$$\matrix{\displaystyle{{1\over{\sqrt{1+\eps}}}}~\N{x}_{G_\th} 
\leq~\N{x}_{E_\th} \leq~\N{x}_{G_\th}~~~\forall x~\in
X\cr
{} \cr
\displaystyle{{1\over{\sqrt{1+\eps}}}}~\N{x}_{H_\th} 
\leq~\N{x}_{F_\th} \leq~\N{x}_{H_\th}~~~\forall x~\in
Y } \eqno (4.2) $$
\par
Let us consider $ x^\ast \in (G_\th)^\ast $  and  $ y^\ast
\in (H_\th)^\ast $  of norm one. Dy the Hahn-Banach
theorem  and (4.2), we can find extensions 
$ \tilde x^\ast \in (E_\th)^\ast $ and 
$ \tilde y^\ast \in (F_\th)^\ast $ with 
$ \tilde x^\ast_{\vert X} =  x^\ast_{\vert X} $ and 
$ \tilde y^\ast_{\vert Y} =  y^\ast_{\vert Y} $ such that
$ \N{ \tilde x^\ast}_{(E_\th)^\ast} \leq~ \sqrt{1 + \eps
}$ ,$ \N{ \tilde y^\ast}_{(F_\th)^\ast} \leq~ \sqrt{1 +
\eps }$. Hence 
$$
\abs{\sum_{k=1}^m ~x^\ast (x_k) y^\ast (y_k) } \leq~
\N{ \tilde x^\ast}_{(E_\th)^\ast}
\N{ \tilde y^\ast}_{(F_\th)^\ast}
\N{T}_{\wco{E_\th}{F_\th}}
\leq (1 + \eps )\N{T}_{\wco{E_\th}{F_\th}}.$$
As $ x^\ast $ and $ y^\ast $ are arbitrary elements on
the unit spheres of $ G_\th^\ast $ and $ H_\th^\ast $ 
 respectively 
$$ \N{T}_{\wco{G_\th}{H_\th}}
\leq (1 + \eps )\N{T}_{\wco{E_\th}{F_\th}}. \eqno
(4.3) $$
\par
Note that if $ Z $ is a closed subspace of a type 2
Banach space $ W $, then the quotient  $ W/Z $ is also 
of type 2, and $ \widetilde T_2( W/Z )
\leq~\widetilde T_2( W ) $. From this we deduce that $
\widetilde T_2( G_k^\ast) \leq~\widetilde T_2( E_k^\ast) $
and $ \widetilde T_2( H_k^\ast) \leq~\widetilde T_2(
F_k^\ast) $ for $ k= 0,1 $.  \par
Using corollary 3.5 we find a constant $ c $ depending
only on $ ~E_0^{}, ~E_1^{}, ~F_0^{}
\hbox{ and } ~F_1^{} $ such that 
$$
\N{T}_{
[\wco{G_0}{H_0}, \wco{G_1}{H_1}]_\th^{}} \leq~
c \N{T}_{\wco{G_\th}{H_\th}}. \eqno  ( 4.4) $$
\par
Since , for $ k= 0,1 $, the canonical embedding $ j :G_k
\wcheck\otimes H_k \rightarrow E_k \wcheck\otimes F_k $
has  a norm smaller than one,  we obtain by interpolation
$$
\N{T}_{
[\wco{E_0}{F_0}, \wco{E_1}{F_1}]_\th^{}}  \leq~
\N{T}_{
[\wco{G_0}{H_0}, \wco{G_1}{H_1}]_\th^{}} . \eqno  ( 4.5)
$$  \par
Since $ \eps $ is arbitrary. Combining (4.3), (4.4) and
(4.5) we obtain the second inequality of (4.1).
\par
In order to prove the first inequality in (4.1),
consider $ T = \sum_{k=1}^m~x_k\otimes y_k \in E 
\otimes F $ and $ \eps > 0 $ as before, and define $ X $
and $ Y $ in the same way. Using lemma 4.1, we can find a
finite-dimensional subspace of $ E \otimes F $
containing 
$ X \otimes Y $, which we can assume to be of the
form  $ G \otimes H $, such that 
$$
\N{T}_{[(G \otimes H)_0 ,( G\otimes H)_1]_\th^{}} 
\leq~ (1 + \eps )
\N{T}_{
[\wco{E_0}{F_0}, \wco{E_1}{F_1}]_\th^{}}. \eqno  ( 4.6) $$ 
But for $ k= 0,1 $, we have 
$(G \otimes H)_k = ( G \otimes H , \N{.}_{\wco{E_0}{F_0}}
) = G_k \wcheck\otimes H_k $. So using the easy part of
corollary 3.5, we obtain 
 $$
\N{T}_{\wco{G_\th}{H_\th}}
\leq
\N{T}_{[(G \otimes H)_0 ,( G\otimes H)_1]_\th^{}}
 \eqno  (4.7) $$ 
Finally, we have trivially
$\N{T}_{\wco{E_\th}{F_\th}} \leq~
\N{T}_{\wco{G_\th}{H_\th}} $. combining this with 
(4.7) and (4.6), the first inequality of (4.1) follows.
This achieves the proof of the theorem.
\par
\bg
In what follows, we will need the following known lemma:
\par
\proclaim Lemma 4.3. 
Let $ V $ and $ W $ be finite-dimensional subspaces  
of $ X $ and $ Y $ respectively. Assume that one of the
following conditions holds \hfill\break
~~~1.~$ X $   and $ Y $ are type 2
spaces.\hfill\break
~~~2.~$ X $  is a type 2 space and $~W $ is a
sublattice of the 2-convex Banach lattice $ Y$.
\hfill\break
~~~3.~$ V $ and $W $ are sublattices
of the 2-convex Banach lattices $~X$ and $~Y$
respectively.\hfill\break
Then, there exists a constant $ c_1 $
depending only on $ X $ and $ Y $ such that
$$ 
\forall~T \in ~V \otimes W ~~~~~
\N{ T }_{X \hot Y} \leq~\N{ T }_{V \hot W} 
\leq~ c_1~\N{ T }_{X \hot Y}.$$ 
\par
\bg
P{\sevenrm ROOF}: 
The first case can be found explicitly in [{\bf M}].
In the other cases, note that if $ u \in {\cal
L}(V , W^\ast )$,  there exist a finite-dimensional
Hilbert space $ H $, and two operators $ A	 : V
\rightarrow H $ and $ B : W \rightarrow H $, such that 
$ u = {}^tB.A $ and $\N{A} = \N{B} \leq c_2 ( \N{u}
)^{1/2} $ with $ c_2 $ a constant depending only on 
type 2 ( or 2-convexity ) constants of $ X $ and $ Y
$.  
Using Maurey's theorem [{\bf M}] in the  type 2 case, or 
well-known facts in the 2-convex Banach lattice case
(see for instance [{\bf Pi4},Proposition 1.1]), we can
find extensions $\widetilde A	 : X
\rightarrow H $ and $\widetilde B : Y \rightarrow H $ of 
the operators $ A $ and $ B $ respectively, such that $
\N{\widetilde A}  \leq c_3 \N{A} $ and $ \N{\widetilde B}
\leq c_3 \N{B} $  
 with $ c_3 $ depending only on type 2 ( or 2-convexity)
constants of $ X $ and $ Y $. If $ \tilde u =
{}^t\widetilde B. \widetilde A \in {\cal L }(X,Y^\ast)$, 
we have $\N{\tilde u } \leq~c_1 \N{u} $ and 
$$ \forall~x \in V~  \forall~y \in W ~~~\langle \tilde
u(x) , y \rangle =\langle u(x) , y \rangle.
 $$
By duality, this implies the difficult inequality of the
lemma.
\par
\bg
\proclaim Theorem 4.4. 
Let $ ( E_0 , E_1 ) $ and $ ( F_0 , F_1 ) $ be two
interpolation pairs. Assume that we are in one of the
following cases
$$\eqalign {~~~1&.~E_0, ~E_1, ~F_0
\hbox{ and } ~F_1~\hbox{ are type 2 spaces }.\cr
~~~2&.~E_0, ~E_1~\hbox{ are type 2 spaces 
and } ~F_0, ~F_1~\hbox{ are 2-convex Banach
lattices}.\cr
~~~3&.~E_0, ~E_1, ~F_0
\hbox{ and } ~F_1~\hbox{ are 2-convex Banach
lattices}. } $$
Then $ ( E_0 \hot F_0 ,E_1 \hot F_1 )$ 
 is an interpolation pair, and for $ 0 < \th < 1 $ we
have
$$
\left[E_0 \hot F_0 ,E_1 \hot F_1\right
]_\th = 
\left[E_0,E_1\right]_\th \hot \left[F_0,F_1\right
]_\th. $$
\par
\bg 
P{\sevenrm ROOF}: 
Let $ X $ ( resp. $ Y $ ) be a type 2 Banach space or a 
 2-convex Banach lattice.
Then using results from [{\bf Pi1},Chapters 3 and 8], we
know that every  operator from $ X $ into $ Y^\ast $
factors through a Hilbert space and consequently, every
operator  from $ X $ into $ Y^\ast $ is approximable. 
This implies that the canonical map 
$ J : X\hot Y \rightarrow X\wcheck\otimes Y $  is
one to one.
\par
Using this remark,  we see  that in each case considered
in the theorem, the spaces $ E_0 \hot F_0 $ and 
 $ E_1 \hot F_1 $ embed continuously in 
$ ( E_0 + E_1 ) \wcheck\otimes ( F_0 + F_1 ) $. Hence 
 $ (E_0 \hot F_0 , E_1 \hot F_1) $ is an
interpolation pair.
\par
With the same convention of notation  as in the begining
of the proof of theorem 4.3, the space $ E \otimes F $ is
clearly  a dense subspace of $ E_\th \hot F_\th $
for every $ \th  \in [0,1] $, and so it is  also  a 
dense   subspace of $ [ E_0 \hot F_0 , E_1
\hot F_1 ]_\th $. 
\par
In order to prove the theorem, it is sufficient to find a
constant $ c $ such that 
$$
\forall~T \in~ E \otimes F~~~
 c^{-1}~\N{T}_{E_\th \hot F_\th} \leq~\N{T}_{
[E_0 \hot F_0 , E_1 \hot F_1]_\th} \leq~
 \N{T}_{E_\th \hot F_\th}. \eqno (4.8) $$
We give only a detailed proof of this inequality  in
the first case. The other cases can be treated similarly
using the corresponding  parts from lemmas 4.1 and 4.3.
\par
In the sequel, we assume that  the hypotheses
of the first case are satisfied.
\par
Let $ T = \sum_{k=1}^m~x_k\otimes y_k \in E \otimes F $ 
and $ \eps > 0 $. Define $ X = \hbox{ span }~\{ x_k
~:~k\leq m \} $ and $ Y = \hbox{ span }~\{ y_k ~:~k\leq
m \} $.  By lemma 4.1, we can find a finite  dimensional 
subspace  of $  E \otimes F $ containing  $ X \otimes Y
$,  which  we can assume to be of the form $ G \otimes H
$, such that
$$
\N{T}_{[(G \otimes H)_0 ,( G\otimes H)_1]_\th} 
\leq~ (1 + \eps )
\N{T}_{
[E_0 \hot F_0 , E_1 \hot F_1]_\th} . \eqno 
 ( 4.9) $$ 
 By lemma 4.3, the identity map $ j : ( G\otimes H)_k
 \rightarrow  G_k \hot H_k  , k =0,1 $  has a
norm  smaller than $ c_1 $  and by interpolation, we
obtain
$$ \N{T}_{
[G_0 \hot H_0 , G_1 \hot H_1]_\th}  \leq~
c_1~ \N{T}_{[(G \otimes H)_0 ,( G\otimes H)_1]_\th} 
 \eqno  ( 4.10) $$ 
Using corollary 3.5  and noting that 
$ {\cal L}(X , Y^\ast ) = ( X \otimes Y )^\ast $  we
have 
$$ \N{T}_{ G_\th \hot H_\th } \leq~c~
\N{T}_{
[G_0 \hot H_0 , G_1 \hot H_1]_\th}  \eqno 
( 4.11 ) $$
Note that the constants $ c $ and $ c_1 $ depend only
on the type 2 constants of the spaces 
$ E_0, E_1, F_0 $ and $ F_1 $.
\par
Trivially, we also have  
$$
\N{T}_{ E_\th \hot F_\th } \leq 
\N{T}_{ G_\th \hot H_\th }   \eqno ( 4.12 ) $$
Combining (4.9), (4.10), (4.11) and  (4.12),  we
obtain the first  inequality of (4.8) which is the
difficult one.
\par
On the other hand, given $ \eps > 0 $  if $ T \in E
\otimes F $,  we can write $ T = \sum_{k=1}^\infty~ 
\alpha_k x_k \otimes y_k $ with $ \alpha_k > 0~,~\N{
x_k}_{E_\th} = \N{ y_k}_{F_\th} = 1 $ and $
\sum_1^\infty \alpha_k \leq~(1 + \eps ) 
\N{T}_{ E_\th \hot F_\th }$.
For every $ k \geq 1 $, we can find  $ P_k $ ( resp. $
Q_k $) an analytic function on the  stripe $ S $,
continuous on $ \adh S $ with values in $ E_0 + E_1 $ 
(resp. $ F_0 + F_1 $ such that 
$ P_k(\th) = x_k $  ( resp.  $ Q_k(\th) = y_k $ ) 
and 
satisfying that  for $ j = 0,1 $ and for 
 every $ t \in \reel $ we have 
$ \N{ P_k (j + it ) }_{E_j} \leq~1 + \eps $
(resp. $ \N{ Q_k (j + it ) }_{F_j} \leq~1 + \eps $.
\par
Let us consider the following function $ U(z) =
\sum_1^\infty~P_k(z) \otimes Q_k(z) $.  This is an analytic 
$( E_0 + E_1 )\wcheck\otimes( F_0 + F_1 )$-valued 
function on $ S $, continuous on $ \adh S $ such that
  $ U(\th) = T $ and   satisfying:
$$
\N{ U(j + it)}_{ E_j \hot F_j } \leq~(1 + \eps)^2~
\sum_{k=1}^\infty \alpha_k \leq~(1 + \eps)^3~
\N{T}_{ E_\th \hot F_\th } $$
for every $ t \in \reel $ and $ j = 0,1 $.
\par
This shows that
$$
\N{T}_{ [E_0 \hot F_0 , E_1 \hot F_1]_\th} 
\leq ~(1 + \eps)^3~
\N{T}_{ E_\th \hot F_\th } $$
and proves  the  second inequality  in (4.8), since $
\eps $ is arbitrary.
\par
\bg
In order to give the right generalization  of  corollary 
3.6, we need some additional  definitions and notation. 
 If $ \{ x_k \}_{k \geq 1} $ is a sequence  of elements  
in a Banach space $ X $ we define 
$$\eqalign { N_2( \{ x_k \}_{k \geq 1} ) &=
\left( \sum_{k=1}^\infty \N{ x_k }^2 \right)^{1/2}. \cr
M_2( \{ x_k \}_{k \geq 1} ) &= \sup \left\{ \left(
\sum_{k=1}^\infty \abs{ \langle \xi , x_k
\rangle}^2\right)^{1/2}~:~\xi \in X^\ast ,~\N{\xi} \leq~1
\right \}. }$$
\par
Let $ u \in E \otimes F $.  Following Saphar in 
[{\bf Sa}], we define the norm 
$$
d_2(u) = \inf \left\{ M_2( \{ x_k \}_{k \geq 1} ).
 N_2( \{ y_k \}_{k \geq 1}) \right\} $$
where the infimum runs over all representations of 
$ u $ of the form $ u = \sum_1^n x_k \otimes y_k $. 
And we denote by
$ E\mathop{\hot}\limits_{d_2} F $ the  completion of 
$ ( E \otimes F, d_2 ) $.
\par
 If we consider the 
elements of $ E \otimes F $ as operators from 
$ F^\ast $ into $ E $ in the natural way, then  by 
[{\bf Sa}]  
$ E \mathop{\hot}\limits_{d_2} F $  is the closure of 
 $ E \otimes F $  in $ \Pi_2(F^\ast , E ) $ and 
$$\forall~~  u \in E \otimes F ~~~d_2(u) = 
\pi_2( u : F^\ast \rightarrow E). $$ 
It is also shown  that the dual space of 
$ E \mathop{\hot}\limits_{d_2} F $ is isometrically isomorphic to 
 $ \Pi_2( E , F^\ast ) $.
\par
Using the preceding notation, we can state the
following  infinite-dimensional generalization of
corollary 3.6.
\par
\proclaim Theorem 4.5. 
Let $ (X_0,X_1 ) $ and  $ (Y_0,Y_1 ) $ be two
interpolation pairs, and assume that 
$ X_0^\ast, X_1^\ast, Y_0^{}$ and $ Y_1^{} $ are type 2 
spaces. Then for $ \th \in ]0,1[ $ we have
$$
\left[ X_0 \mathop{\hot}\limits_{d_2} Y_0 ,X_1
 \mathop{\hot}\limits_{d_2}
Y_1 \right]_\th = 
[X_0,X_1]_\th \mathop{\hot}\limits_{d_2} [Y_0,Y_1]_\th.  $$
\par
\bg
{\it Remark.}
 A Banach lattice version of this theorem also holds, we
leave the detailed formulation as an exercise for the
interested reader.
\par
P{\sevenrm ROOF}:(sketch). 
Since $  X_k \mathop{\hot}\limits_{d_2} Y_k $ can be considered
as a closed  subspace of $ \Pi_2(Y_k^\ast , X_k^{} ) $
$(\hbox{ for }k = 0,1 )$, then clearly both $  X_0
\mathop{\hot}\limits_{d_2} Y_0 $ and $  X_1 
\mathop{\hot}\limits_{d_2} Y_1 $
embed in  ${\cal L}(Y_0^\ast \cap Y_1^\ast , X_0^{} +
X_1^{} ) $ and $ (X_0 \mathop{\hot}\limits_{d_2}
Y_0, X_1 \mathop{\hot}\limits_{d_2} Y_1) $ is an interpolation
pair.
\par
Clearly, Corollary 3.6 can be formulated in the following 
way: There exists a constant $ c > 0 $ ( depending only
on the type 2 constants of $ X_0^\ast, X_1^\ast, Y_0^{} $
and $ Y_1^{} $ ) such that if $ M $ and $ L $ are 
finite-dimensional subspaces of $ X_0 \cap X_1 $  and 
$ Y_0 \cap Y_1 $ respectively, then for all $ T \in M
\otimes L $ 
$$
{1\over c}\N{T}_{ M_\th \mathop{\hot}\limits_{d_2} L_\th }\leq~
\N{T}_{[ M_0 \mathop{\hot}\limits_{d_2} L_0 , M_1 
\mathop{\hot}\limits_{d_2}
L_1]_\th }\leq~\N{T}_{ M_\th \mathop{\hot}\limits_{d_2} L_\th } $$
\par
Now, the proof of the theorem goes on exactly as the proof
of  theorem 4.4, using the following lemma 4.6 instead of
lemma 4.3.
\par
\proclaim Lemma 4.6. 
Let  $ M $ and $ L $ be finite-dimensional subspaces of 
$ X $ and $ Y $ respectively, and assume that $ Y  $ is
a type 2 space. Then, for all $ T \in M \otimes L $ 
$$
\N{T}_{ X \mathop{\hot}\limits_{d_2} Y} \leq~
\N{T}_{ M \mathop{\hot}\limits_{d_2} L} \leq~ \widetilde T_2(Y)
\N{T}_{ X \mathop{\hot}\limits_{d_2} Y}.$$
\par
This follows from the fact that for every $ T \in 
\Pi_2(M , L^\ast ) $,  there exists $ \widetilde T \in 
\Pi_2(X , Y^\ast ) $ such that  
$$
\pi_2( \widetilde T : X \rightarrow Y^\ast ) \leq~
\widetilde T_2(Y)~\pi_2( \widetilde T : M \rightarrow
L^\ast ) $$ and
$$
\forall~x \in M~, ~\forall~y \in L ~~~~
\langle T(x),y \rangle = \langle \widetilde T(x),y
\rangle .$$  This,  in turn, follows from Pietch's
factorization theorem, and Maurey's extension theorem
[{\bf M}].
\par
\bg
\noindent{\bf 5. Applications }
\bg
Let us give some consequences of the preceding results.
\par
In the following corollary $ L^p $ will denote 
$ L^p([0,1],dt) $ for $ p \in [1,+\infty[ $ and  the space
$ L_0^\infty([0,1],dt) $ if $ p = +\infty $ (see [{\bf
BL},chapter 5] for the notation).
\par
\proclaim Corollary 5.1.\sn
1. If $ p_0, p_1, q_0 $  and $ q_1$ are
elements of $[1,2] $, and $ \th \in ]0,1[ $, then
$$\left[ L^{p_0}\wcheck\otimes L^{q_0} ,
L^{p_1}\wcheck\otimes L^{q_1} \right]_\th =
 L^{p_\th}\wcheck\otimes L^{q_\th}. $$
\indent 2. If $ p_0, p_1, q_0 $  and $ q_1$ are
elements of $[2,+\infty] $, and $ \th \in ]0,1[ $, then
$$\left[ L^{p_0}\hot L^{q_0} ,
L^{p_1}\hot L^{q_1} \right]_\th =
 L^{p_\th}\hot L^{q_\th}. $$
With $\displaystyle{ {1\over{p_\th}} = {{1-\th}\over{p_0}} +
{\th\over{p_1}} ,{1\over{q_\th}} = {{1-\th}\over{q_0}} +
{\th\over{q_1}}} $.
\par
\bg
Using the results of [{\bf Z}], one can show that for any  
$ p \in [1,2] $ and $ \th \in ]0,1[ $ we have 
$[L^2\hot L^2 , L^2\hot L^1]_\th \not = 
L^2\hot L^p $. ( The proof is based upon the fact
that for every even integer $ q $, there exists a $
\Lambda (q) $ set which is not a $ \Lambda (q + \eps) $ 
set for any $ \eps > 0 $). This shows that our results can
not be significantly improved and that the assumptions of 
theorems 4.2 and 4.4 are essential.
\par
{\it Remark.} 
In our framework for the interpolation of Banach lattices
we excluded the case of $ C(K) $. It is not difficult to
see that we can take for the spaces $ E_0 $ or $ F_0 $ a 
$ C(K)$-space in theorem 4.4,  and  take $ C([0,1]) $ 
instead of $ L_0^\infty([0,1]) $ in corollary 5.1.
\par
We also have a non-commutative version of the
preceding corollary.
\par
Let us recall first the following definition of unitary
ideals. If $ A $ is a compact operator acting on a Hilbert
space, then $ \abs{A} $ denotes the modulous of $ A $,\ie 
$\abs{A} = \sqrt{ A^\ast A } $. And $ s(A) = \{
s_n(A)\}_{n \geq 1} $ denotes the sequence of singular
numbers of $ A $,\ie $ s_n(A) $ is the $n$th  eigenvalue
of $ \abs{A} $ ( where eigenvalues are counted in
non-increasing order, according to their multiplicity).  
Suppose that $ (E,\N{.}_E) $ is a symmetric Banach
sequence  space.  The corresponding unitary ideal $ C_E $ 
is the space
$$ C_E = \{ A~\hbox{  compact }~~: ~s(A) \in E \} $$
with the norm $ \N{A}_{C_E} = \N{s(A)}_E $  for $ A \in
C_E $.
\par
In the case $ E = \ell_p $,  for $ p \in [1,+\infty[ $  
we use the notation $ C_{\ell_p} = C_p $ and 
$\N{A}_p = \left(\sum_1^\infty (s_k(A))^p\right)^{1/p} $  
for $ A \in C_p $. Finally, $ C_\infty $ denotes the
space of all compact operators equipped with the usual
operator norm.
\par
Using the results of [{\bf T}] we know that $ C_p $ is
of  type 2 if $ p\in [2,+\infty[ $. It is also known that 
$[C_{E_0},C_{E_1}]_\th = C_{E_\th} $ if $ E_0 $, $ E_1$ 
are symmetric Banach sequence spaces, and $ E_\th =
[E_0,E_1]_\th $.
\par
Exploiting these facts, we obtain the following
non-commutative analogue of corollary 5.1.
\proclaim Corollary 5.2.\sn
1. If $ p_0, p_1, q_0 $  and $
q_1$ are elements of $]1,2] $, and $ \th \in ]0,1[$, then
$$\left[ C_{p_0}\wcheck\otimes C_{q_0} ,
C_{p_1}\wcheck\otimes C_{q_1} \right]_\th =
 C_{p_\th}\wcheck\otimes C_{q_\th}. $$
\indent 2. If $ p_0, p_1, q_0 $  and $ q_1$ are
elements of $[2,+\infty[ $, and $ \th \in ]0,1[$, then
$$\left[ C_{p_0}\hot C_{q_0} ,
C_{p_1}\hot C_{q_1} \right]_\th =
 C_{p_\th}\hot C_{q_\th}. $$
With $\displaystyle{ {1\over{p_\th}} = {{1-\th}\over{p_0}} +
{\th\over{p_1}} ,{1\over{q_\th}} = {{1-\th}\over{q_0}} +
{\th\over{q_1}}} $.
\par
\bg
As an application of theorem 4.5, we will show the
following result, which yields a non-commutative analogue 
of the fact that if $ 1 < p \leq 2 $, then $
\Pi_2(\ell_p,\ell_p) $ is super-reflexive [{\bf Pi3}].
\par
\proclaim Theorem 5.3. 
Let $ E $ be a type 2 symmetric Banach sequence space.
 Then, there exists an interpolation pair $ ({\cal H},X)
$,  where ${\cal H} $  is a Hilbert space, and $ \th \in
]0,1[ $ such that 
$$
\Pi_2 ( C_{E^\ast},C_{E^\ast} ) = [{\cal H},X]_\th. $$
\par
\bg
P{\sevenrm ROOF}: 
Since $ E $ has type 2,  $ E $ is 2-convex and
$q$-concave for some $ q \in [2,+\infty[ $. Using theorem 
2.2 and remark 2.5 of [{\bf Pi5}], we can find 
$ E_0 $ a  2-convex symmetric Banach sequence space and 
$ \th \in ]0,1[ $ such that $ E = [\ell_2,E_0]_\th $.  
If, moreover, we use the reiteration theorem of
interpolation and modify the value of $ \th \in ]0,1[ $,
if necessary, then we can assume that $ E_0 $ is 
$K$-convex, hence that $ E_0 $ is of type 2.
\par
It follows that $ C_E = [ C_2 , C_{E_0}]_\th $
  and we know by the results of [{\bf GT}] that 
$  C_{E_0} $ is of type 2. Applying theorem 4.5 to the
pairs $ ( C_2 , C_{E_0^\ast} ) $ and $ ( C_2 , C_{E_0^{}})
$  we obtain
$$
C_{E^\ast} \mathop{\hot}\limits_{d_2} C_E =
\left[ C_2 \mathop{\hot}\limits_{d_2} C_2 , C_{E_0^\ast}
\mathop{\hot}\limits_{d_2} C_{E_0} \right]_\th. $$
By  duality, we obtain 
$$
\Pi_2 ( C_{E^\ast},C_{E^\ast} ) = 
\left[\Pi_2 ( C_2,C_2 ) ,\Pi_2 (
C_{E_0^\ast},C_{E_0^\ast} ) \right]_\th. $$
Which implies the result since $ \Pi_2 ( C_2,C_2 ) $ is a 
Hilbert space.
\par
\proclaim Corollary 5.4. 
Let $ E $ be a type 2 symmetric Banach sequence space.
Then
$\Pi_2 ( C_{E^\ast},C_{E^\ast} ) $  
is super-reflexive. In particular, if $ 1 < p \leq 2 $
then $\Pi_2 ( C_p,C_p ) $ is super-reflexive.
\par
\bg
Note that, by the results of [{\bf L}], if $ p > 2 $ then 
$\Pi_2 ( C_p,C_p ) $ is not even $ B$-convex.
\par
Argueing exactly as in corollary 3.6 and using theorem
4.2, we obtain the following corollary:
\par
\proclaim Corollary 5.5. 
Let $( E_0 , E_1 ) $ and  $( F_0 , F_1 ) $ be two
interpolation pairs, and assume that $ E_0^\ast ,
E_1^\ast $ are type 2 space or  that $  E_0 , E_1 $ are 
2-concave Banach lattices. Then, for every $ p_0, p_1 \in 
[1,2] $ and $ \th \in ]0,1[ $ the following holds
$$
\left[\Pi_{p_0}( E_0 , F_0 ) , \Pi_{p_1}( E_0 , F_1 )
\right]_\th \subset \Pi_p([E_0 , E_1 ]_\th,[F_0,F_1 ]_\th
). $$
With $\displaystyle{ {1\over{p_\th}} = {{1-\th}\over{p_0}} +
{\th\over{p_1}} }$.
\par
\bg
Now  we  give  another application of our results. If 
$ X $ is a complex Banach space and $ 0 < p < +\infty $, 
then $ \widetilde H^p(X) $ denotes the closure in $ L^p(
\tore;X) $ of the set of analytic $X$-valued polynomials:
 $$
\left\{ \sum_{k=0}^n~z^k~x_k ~: ~ n \in \nat ,~z \in \tore
 \hbox{ and } x_k  \in X \right \}. $$
The space $ \widetilde H^p(X) $ equipped with the norm 
(quasi-norm if $ p < 1 $ ) induced by $ L^p( \tore ; X ) $ 
  is a Banach space ( quasi-Banach space).
\par
Let us formulate the following definition:
\par
\bg
D{\sevenrm EFINITION}. 
An interpolation pair $ (X_0 , X_1) $ of complex Banach
spaces will be called a Hardy-interpolation pair, if and
only if, for $ p_0, p_1 \in [1,+\infty[ $ and $ \th \in 
]0,1[ $
$$
\widetilde H^{p_\th}([X_0,X_1]_\th ) \subset 
\left[ \widetilde H^{p_0}(X_0) , \widetilde H^{p_1}(X_1)
\right]_\th, $$
with $\displaystyle{ {1\over{p_\th}} = {{1-\th}\over{p_0}} +
{\th\over{p_1}}} $.
\par
As it is noted in [{\bf BX}], the converse inequality 
always holds true with no assumptions on the spaces, and 
actually, if $ (X_0 , X_1) $  is  a Hardy-interpolation
pair then for all $ p_0, p_1 \in [1,+\infty[ $ and $ \th
\in  ]0,1[ $
$$
\widetilde H^{p_\th}([X_0,X_1]_\th ) =
\left[ \widetilde H^{p_0}(X_0) , \widetilde H^{p_1}(X_1)
\right]_\th, $$
with $\displaystyle{ {1\over{p_\th}} = {{1-\th}\over{p_0}} +
{\th\over{p_1}}} $.
\par
\bg
If $ X $ is a complex banach space, then $ (X , X ) $ is a 
 Hardy-interpolation pair. This was noted by G. Pisier,
see [{\bf BX}]. It is also shown in [{\bf BX}] that if 
$ X_0 $ and $ X_1 $ are UMD spaces then $ (X_0 , X_1 ) $ 
is a  Hardy-interpolation pair, whereas $ ( L^1(\tore),
c_0(\ent) ) $ is not a  Hardy-interpolation pair.
\par
We think that it would be interesting  to study the 
class of  Hardy-interpolation pairs, and one step in this 
direction is to study the stability properties of this
class.
\par
It is not difficult to see that if $ (X_0 , X_1 ) $ 
is a  Hardy-interpolation pair then, for every $ 
r_0, r_1 \in [1,+\infty[ $ the pair 
$ ( L^{r_0}(X_0) , L^{r_1}(X_1) ) $ is a 
 Hardy-interpolation pair where $  L^{r}(X) = L^r(\Omega,
{\cal A}, \pe; X) $.
\par
We also have the following:
\par
\proclaim Theorem 5.6. 
Let $ X_0,  X_1, Y_0 $ and $ Y_1 $ be complex Banach
spaces of type 2. Assume that $ (X_0 , X_1) $ and 
$ (Y_0 , Y_1) $ are two  Hardy-interpolation pairs. 
Then  $(X_0\hot Y_0 , X_1\hot Y_1 ) $ is
also a  Hardy-interpolation pair.
\par
\bg
P{\sevenrm ROOF}: 
Let us consider $ p_0, p_1 \in[1,+\infty[
$, then the pairs $\left( \widetilde H^{2p_0}(X_0)
, \widetilde H^{2p_1}(X_1)\right) $ and 
$\left( \widetilde H^{2p_0}(Y_0) , \widetilde
H^{2p_1}(Y_1)\right) $ are two interpolation pairs of
type 2 spaces. Using theorem 4.4 and the hypothesis we
obtain, for every $ 0 < \th <1 $,
\bg 
$$  \left[\widetilde H^{2p_0}(X_0)\hot 
\widetilde H^{2p_0}(Y_0)  , \widetilde
H^{2p_1}(X_1) \hot \widetilde
H^{2p_1}(Y_1)\right]_\th =~~~~~~~~~~~~~~~~~~~~~~~~~~~~~~$$
$$\eqalign {~~~~~~~~~~~~~ =&\left[\widetilde H^{2p_0}(X_0)
, \widetilde H^{2p_1}(X_1)\right]_\th \hot 
\left[ \widetilde H^{2p_0}(Y_0) , \widetilde
H^{2p_1}(Y_1)\right]_\th  \cr
=&\widetilde H^{2p_\th}(X_\th)\hot 
\widetilde H^{2p_\th}(Y_\th).} $$
\par
Let $ \widetilde Q : \widetilde H^{2r}(V)\hot 
\widetilde H^{2r}(W) \rightarrow \widetilde H^r( V\hot W)
$ be the natural norm one operator considered by  Pisier
in [{\bf Pi2}] : $$
\widetilde Q \left ( \sum_{k=1}^n f_k \otimes g_k
\right)(z) = \sum_{k=1}^n f_k(z) \otimes g_k(z).$$
By interpolation  
$$
\widetilde Q\left( \widetilde H^{2p_\th}(X_\th)\hot 
\widetilde H^{2p_\th}(Y_\th) \right) \subset 
\left[ \widetilde H^{p_0}(X_0 \hot Y_0 ), 
\widetilde H^{p_1}(X_1 \hot Y_1 )\right ]_\th.
$$
By theorem 3.1 of [{\bf Pi2}], since $ X_\th $ and 
$ Y_\th $ have type 2, we obtain
$$
\widetilde Q\left( \widetilde H^{2p_\th}(X_\th)\hot 
\widetilde H^{2p_\th}(Y_\th) \right) = 
\widetilde H^{p_\th}(X_\th \hot Y_\th ). $$
Hence,
$$ 
\widetilde H^{p_\th}(X_\th \hot Y_\th ) \subset 
\left[ \widetilde H^{p_0}(X_0 \hot Y_0 ), 
\widetilde H^{p_1}(X_1 \hot Y_1 )\right ]_\th.$$
And the result follows by a further use of theorem 4.4.
\bg
{\it Acknowledgement :}
\par
The author thanks Professor Gilles Pisier for many
stimulating discussions during the preparation of this
work.
\bg
\centerline{\bf R{\sevenbf EFERENCES}}
\bg
\refl AB&M.H. Aizenshtein and Yu.A. Brudnyi& Computable
interpolation functions.&Investigation on the theory of
functions of many real variables.& Jaroslavl Univ.
(1986), 11--35. (Russian).
\par
\refl Ai&M.H. Aizenshtein & duality of the interpolation
functions.& ibid & 3--10 (Russian).
\par
\refl BL&J. Berg and J. L\"ofstr\"om && Interpolation
spaces, an introduction.&Springer-Verlag, Berlin and
New-York,(1976).
\par
\refl BX&O. Blasco and Q. Xu & Interpolation between
vector valued Hardy spaces.&{ } & To Appear.
\par
\refl C&A.P. Calder\'on & intermediate
spaces  and interpolation, the complex method.&Studia
Math.& 24 (1964),113--190.
\par
\refl CCRSW1&R. Coifman, M. Cwikel, R. Rochberg, Y.
Sagher and G. Weiss & The complex method for
interpolation of operators acting on families of Banach
spaces.&Lecture Notes in Mathematics No 779&
Springer,(1980).
\par
\refl CCRSW2&R. Coifman, M. Cwikel, R. Rochberg, Y.
Sagher and G. Weiss & A theory of complex interpolation
for families of Banach spaces.& Advances in Math.&vol.43,
No 3 (1982),203--229.
\par
\refl CS&R. Coifman and S. Semmes & Interpolation of
Banach spaces, Perron process and Yang-Mills.&
{}&Preprint.
\par
\refl GT&D.J.H. Garling and N. Tomczak-Jaegermann& The
cotype and uniform convexity of unitary ideals.&Is.
J. Math.& vol.45 No 2-3 (1983),175--197.
\par
\refl H &H. Helson&& Lectures on invariant
subspaces.&Academic Press, (1964).
\par
\refl He&E. Hernandez& Intermediate spaces and the
complex method of interpolation for families of Banach
spaces.&Ann. Scu. Norm. Sup. Pisa S\'erie IV &vol.13 No 2
(1986).
\par
\refl K &O. Kouba& Sur l'interpolation des produits
tensoriels projectifs ou injectifs d'espaces de Banach&
C.R.Acad.Sci. Paris S\'erie I& t.309, No
10,(1989),683--686.
\par
\refl K\"o&G. K\"othe&& Topological vector spaces II.&
Springer-Verlag, Berlin and New-York,(1979).
\par
\refl KP& S. Kaijser and J.W. Pelletier& Interpolation
functors and duality.& Lecture Notes in Mathematics No
1208&Springer-verlag,(1986).
\par
\refl Kw&S. Kwapie$\acute{\hbox{n}}$ & Isomorphic
characterization of inner product spaces by orthogonal
series with vector coefficients.&Studia Math.&
44,(1972),583--595.
\par
\refl L&P.K. Lin&B-convexity of the space of 2-summing
operators.&Is. J. Math.&vol. 37,(1980),139-150.
\par
\refl LT1&J. Lindenstrauss and L. Tzafriri&& Classical
Banach spaces I.& Springer-Verlag, Berlin and
New-York,(1977).
\par 
\refl LT2&J. Lindenstrauss and L. Tzafriri&& Classical
Banach spaces II.& Springer-Verlag, Berlin and
New-York,(1979).
\par                                              
\refl M& B. Maurey& Un th\'eor\`eme de prolongement.&
C.R. Acad. Sci. Paris, S\'erie I.&t.279,(1974),329--332.
\par
\refl P&A. Pietsch&& Operator ideals.&V.E.B. Berlin and
North-Holland,\  Amsterdam (1978).
\par
\refl Pi1&G. Pisier&& Factorization of linear operators 
and geometry of Banach spaces.&C.B.M.S. No 60, Amer.
Math. Soc. Providence, second printing (1987).
\par
\refl Pi2&G. Pisier& Factorization of operator valued
analytic functions.& Advances in Math.&  93 No 1, (1992) 61--125.
\par
\refl Pi3&G. Pisier& A remark on $ \Pi_2(\ell^p,\ell^p) $.&
Matematische Nachrichten.& 148 (1990), 243-245.
\par
\refl Pi4&G. Pisier&& Topics on Grothendieck's theorem.
 Proceedings of the international conference on operator
algebras, ideals, and their applications in theoretical
Physics.&Teulner
Verlagsgesellschaft, Leipzig, (1978), 44--57.
 \par
\refl Pi5&G. Pisier&Some applications of the complex
interpolation method to Banach lattices.&J. Analyse
Math.&35, (1979),264--281.
\par
\refl Sa&P. Saphar& Produits tensoriels d'espaces de\  
Banach\  et classes d'appli- cations lin\'eaires.&Studia
Math.&38, (1970),71--100.
\par
\refl T&N. Tomczak-Jaegermann&The moduli of smoothness and 
convexity and the Rademacher averages of the trace
classes $ S_p~( 1\leq p <\infty )$.&
Studia Math.&50, (1974), 163--182.
\par
\refl Z&M. Zafran& Interpolation of multiplier spaces.& 
Amer. J. Math.&vol.105,No 6, (1983),1405--1416.
\par

\end